\documentclass[a4paper,10pt]{amsart}

\title{Trace ideals for Fourier integral operators with
non-smooth symbols II}
\usepackage{cite}
\usepackage{amssymb}
\usepackage{latexsym}
\usepackage{amsmath}
\usepackage{mathrsfs}
\usepackage[X2,T1]{fontenc}
\usepackage{amsthm}
\usepackage{calc}                         

\newcommand{\scal}[2]{\langle #1,#2\rangle}

\newcommand{\rr}[1]{\mathbf R^{#1}}

\newcommand{\nm}[2]{\Vert #1\Vert _{#2}}
\newcommand{\nmm}[1]{\Vert #1\Vert }
\newcommand{\abp}[1]{\vert #1\vert}

\newcommand{\op}{\operatorname{Op}}
\newcommand{\opp}{\operatorname{\mathsf{Op}}}

\newcommand{\fy}{\varphi}

\newcommand{\cdo}{\, \cdot \, }

\newcommand{\supp}{\operatorname{supp}}

\newcommand{\vrum}{\vspace{0.1cm}}

\numberwithin{equation}{section}          
\newtheorem{thm}{Theorem}
\numberwithin{thm}{section}
\newtheorem*{tom}{\rubrik}
\newcommand{\rubrik}{}
\newtheorem{prop}[thm]{Proposition}

\newtheorem{lemma}[thm]{Lemma}

\theoremstyle{definition}

\newtheorem{defn}[thm]{Definition}

\theoremstyle{remark}

\newtheorem{rem}[thm]{Remark}              

\author{Francesco Concetti}

\address{Department of Mathematics,
Turin University, Italy}

\email{francesco.concetti@unito.it}

\author{Gianluca Garello}

\address{Department of Mathematics,
Turin University, Italy}

\email{gianluca.garello@unito.it}

\author{Joachim Toft}

\address{Department of Mathematics and Systems Engineering,
V{\"a}xj{\"o} University, Sweden}

\email{joachim.toft@vxu.se}

\begin{document}

\begin{abstract}
We consider Fourier integral operators with symbols in modulation
spaces and non-smooth phase functions whose second orders of
derivatives belong to certain types of modulation space. We establish
continuity and Schatten-von Neumann properties of such operators when
acting on modulation spaces.
\end{abstract}

\maketitle

\section{Introduction}\label{sec0}

\par

In \cite{Bu1}, A. Boulkhemair considers a certain class of Fourier
integral operators were the corresponding symbols are defined without
any explicit regularity assumptions and with only small
regularity assumptions on the phase functions. The symbol class here
is, in the present paper, denoted by $M^{\infty ,1}$ and
contains $S^0_{0,0}$, the set of smooth functions which are bounded
together with all their derivatives. In the time-frequency
community, $M^{\infty ,1}$ is known as a modulation space. (See
e.{\,}g. \cite {Gc2, Fe4, FG1} or the definition below.) Boulkhemair
then proves that such operators are uniquely extendable to continuous
operators on $L^2$. In particular it follows that pseudo-differential
operators with symbols in $M^{\infty ,1}$ are $L^2$-continuous, which
was proved by J. Sj{\"o}strand in \cite {Sj1}, where it seems that
$M^{\infty ,1}$ was used for the first time in this context.

\par

More recent contribution to the theory of Fourier integral operators
with non-smooth symbols are presented in \cite{RuS1,RuS2,RuS3}. For
example, in \cite{RuS2}, Ruzhansky and Sugimoto investigate, among
others, $L^2$ estimates for Fourier integral operators with symbols in
local Sobolev-Kato spaces, and with less regularity assumptions on the
phase functions comparing to \cite{Bu1}.

\par

In this paper we consider Fourier integral operators were the symbol
classes are given by $M^{p,q}$ where $p,q\in [1,\infty]$, and with phase
functions satisfying similar conditions as in \cite{Bu1}. We
discuss continuity of such operators when acting on modulation
spaces, and prove Schatten-von Neumann properties when acting on
$L^2$.

\par

In order to be more specific we recall some definitions. Assume that
$p,q\in [1,\infty ]$ and that $\omega \in \mathscr P(\rr {2n})$ (see
Section \ref{sec1} for the definition of $\mathscr P(\rr n)$). Then
the \emph{modulation space} $M^{p,q}_{(\omega )}(\rr
n)$ is the set of all $f\in \mathscr S'(\rr n)$ such that
\begin{equation}\label{modnorm}
\nm {f}{M^{p,q}_{(\omega )}} \equiv \Big( \int \Big( \int|\mathscr F(f
\tau_x\chi)(\xi)\omega (x,\xi )|^p\, dx\Big)^{q/p}\, d\xi
\Big)^{1/q}<\infty
\end{equation}
(with obvious modification when $p=\infty$ or
$q=\infty$). Here $\tau _x$ is the translation operator $\tau _x\chi
(y)=\chi (y-x)$, $\mathscr F$ is the Fourier transform on $\mathscr
S'(\rr n)$ which is given by
$$
\mathscr Ff(\xi ) = \widehat f(\xi )\equiv (2\pi )^{-n/2}\int
f(x)e^{-i\scal x\xi}\, dx
$$
when $f\in \mathscr S(\rr n)$, and $\chi \in \mathscr S(\rr
n)\setminus 0$ is called a \emph{window function} which is kept
fixed. For conveniency we set $M^{p,q}=M^{p,q}_{(\omega )}$ when
$\omega =1$.

\par

During the last twenty years, modulation spaces have been an active
fields of research (see e.{\,}g. \cite {Gc2, Fe4, Fe5, FG1, FG2,
PT1, Te, To7}). They are rather similar to Besov spaces (see
\cite{BaC,ST,To7} for sharp embeddings) and it has appeared that they are
useful to have in background in time-frequency analysis and to some
extent also in pseudo-differential calculus.

\par

Next we discuss the definition of Fourier integral operators. For
conveniency we restricts ourself to operators which belong to
$\mathscr L(\mathscr S(\rr n),\mathscr S'(\rr n))$. Here we let
$\mathscr L(V_1,V_2)$ denote the set of all linear and continuous
operators from $V_1$ to $V_2$, when $V_1$ and $V_2$ are topological
vector spaces. For any appropriate $a \in \mathscr S'(\rr {2n+m})$ (the
\emph{symbol}) and real-valued $\fy \in C(\rr {2n+m})$ (the \emph{phase
function}), the Fourier integral operator $\op _\fy (a)$ is defined by
the formula
\begin{equation}\label{fintop}
\op _\fy (a)f(x) =(2\pi )^{-n}\iint a(x,y,\xi )f(y)e^{i\fy (x,y,\xi
)}\, dyd\xi ,
\end{equation}
when $f\in \mathscr S(\rr n)$. Here the integrals should be
interpreted in distribution sense if necessary. By letting $m=n$, and
choosing symbols and phase functions in appropriate ways, it
follows that the pseudo-differential operator
$$
\op (a)f(x) =(2\pi )^{-n}\iint a(x,y,\xi )f(y)e^{i\scal{x-y}\xi}
\, dyd\xi 
$$
is a special case of Fourier integral operators. Furthermore, if
$t\in \mathbf R$ is fixed, and $a$ is an appropriate function or
distribution on $\rr {2n}$ instead of $\rr {3n}$, then the definition
of the latter pseudo-differential operators cover the definition of
pseudo-differential operators of the form
\begin{equation}\label{e0.5}
a_t(x,D)f(x) = (2\pi )^{-n}\iint a((1-t)x+ty,\xi
)f(y)e^{i\scal{x-y}\xi}\, dyd\xi .
\end{equation}

\par

On the other hand, in the framework of harmonic analysis it follows
that the map $a\mapsto a_t(x,D)$ from $\mathscr S(\rr {2n})$ to
$\mathscr L(\mathscr S(\rr n),\mathscr S'(\rr n))$ is uniquely
extendable to a bijection from $\mathscr S'(\rr {2n})$ to $\mathscr
L(\mathscr S(\rr n),\mathscr S'(\rr n))$. Consequently, any Fourier
integral operator is equal to a pseudo-differential operator
of the form $a_t(x,D)$.

\par

In the litterature it is usually assumed that $a$ and $\fy$ in
\eqref{fintop} are smooth functions. For example, if $a\in \mathscr
S(\rr {2n+m})$ and $\fy \in
C^\infty (\rr {2n+m})$ satisfy $\partial ^\alpha \fy\in S^0_{0,0}(\rr
{2n+m})$ when $|\alpha |=N$ for some integer $N\ge 0$, then it is easily
seen that $\op _\fy (a)$ is continuous on $\mathscr S(\rr n)$ and
extends to a continuous map from $\mathscr S'(\rr n)$ to $\mathscr
S(\rr n)$. In \cite{AF} it is proved that if $\partial ^\alpha \fy\in
S^0_{0,0}(\rr {2n+m})$ when $|\alpha |=2$ and satisfies
\begin{equation}\label{detphicond}
\left | \det \left ( \begin{matrix}
                     \fy ''_{x,y} & & \fy ''_{x,\xi}
\\[1ex]
                     \fy ''_{y,\xi } &  &\fy ''_{\xi ,\xi}
                     \end{matrix}
\right ) \right | \ge \mathsf d
\end{equation}
for some $\mathsf d >0$, then the definition of $\op _\fy$ extends
uniquely to any $a\in S^0_{0,0}(\rr {2n+m})$, and then $\op _\fy (a)$
is continuous on $L^2(\rr n)$. Next assume that $\fy$ instead
satisfies $\partial ^\alpha \fy\in M^{\infty ,1}(\rr {3n})$ when
$|\alpha |=2$ and that \eqref{detphicond} holds for some $\mathsf d
0$. This implies that the condition on $\fy$ is relaxed since
$S^0_{0,0} \subseteq M^{\infty ,1}$. Then Boulkhemair improves the
result in \cite{AF} by proving that the definition of $\op _\fy$
extends uniquely to any $a\in M^{\infty ,1}(\rr {2n+m})$, and that $\op
_\fy (a)$ is still continuous on $L^2(\rr n)$.

\par

In Section \ref{sec2} we discuss Schatten-von Neumann properties for
Fourier integral operators which are related to those which were
considered by Boulkhemair. More precisely, let $p'\in [1,\infty ]$
denote the conjugate exponent of $p\in [1,\infty ]$. Then we prove
that if $\omega _1,\omega _2\in \mathscr P(\rr {2n})$ and $\omega \in
\mathscr P(\rr {4n})$ are appropriate weight functions, $p,q\in
[1,\infty ]$ are such that $q\le \min (p,p')$ and $a\in
M^{p,q}_{(\omega )}(\rr {2n})$ then $\op _\fy (a)$ belongs
to $\mathscr I_p(\mathscr H_1,\mathscr H_2)$, the set of
Schatten-von Neumann operators of order $p\in [1,\infty]$ from the
Hilbert space $\mathscr H_1=M^2_{(\omega _1)}$ to $\mathscr
H_2=M^2_{(\omega _2)}$. Recall that an operator $T$ from $\mathscr
H_1$ to $\mathscr H_2$ is a Schatten-von
Neumann operator of order $p$ if it is linear and continuous from
$\mathscr H_1$ to $\mathscr H_2$, and satisfies 
$$
\nm T{\mathscr I_p}\equiv \sup \Big ( \sum|(Tf_j,g_j)_{\mathscr H_2}|^p\Big )
^{1/p}<\infty .
$$
Here the supremum should be taken over all orthonormal
sequences $(f_j)$ in $\mathscr H_1$ and $(g_j)$ in $\mathscr H_2$.

\par

Furthermore, assume that $p,q\in [1,\infty ]$ are such that $\le p$
and $p\le p'$, $m=n$ and instead $a(x,y,\xi
)=b(x,\xi )$, for some $b\in M^{p,q}(\rr {2n})$, and that in addition
\begin{equation}\label{detphicond2}
|\det (\fy ''_{y,\xi})|\ge \mathsf d
\end{equation}
holds for some constant $\mathsf d >0$. Then we prove that $\op _\fy
(a)\in \mathscr I_p$. When proving these results we first prove that the
they hold in the case $p=1$. The remaining cases are then consequences
of Boulkhemair's result, interpolation and duality.

\par

\section{Preliminaries}\label{sec1}

\par

In this section we discuss basic properties for modulation
spaces. The proofs are in many cases omitted since they can
be found in \cite {Fe2, Fe3, Fe4, Fe5, FG1, FG2, FG4, Gc2,
To1, To2, To7}.

We start by discussing some notations. The duality between a
topological vector space and its dual
is denoted by $\scal \cdo \cdo $. For admissible $a$ and $b$ in $\mathscr
S'(\rr n)$, we set $(a,b)=\scal a{\overline b}$, and it is obvious that
$(\cdo ,\cdo )$ on $L^2$ is the usual scalar product.

\par

Next assume that $\mathscr B_1$ and $\mathscr B_2$ are topological
spaces. Then $\mathscr B_1\hookrightarrow \mathscr B_2$ means that
$\mathscr B_1$ is continuously embedded in $\mathscr B_2$. In the case
that $\mathscr B_1$ and $\mathscr B_2$ are Banach spaces, $\mathscr
B_1\hookrightarrow \mathscr B_2$ is equivalent to $\mathscr
B_1\subseteq \mathscr B_2$ and $\nm x{\mathscr B_2}\le C\nm x{\mathscr
B_1}$, for some constant $C>0$ which is independent of $x\in \mathscr
B_1$.

\medspace

Next we discuss appropriate conditions for the involved weight
functions. Let $\omega, v\in L^\infty _{loc}(\rr {n})$ be positive
functions. Then $\omega$ is called $v$-\textit{moderate} if
\begin{equation}\label{moderate}
\omega(x+y) \leq C \omega(x) v(y), \quad x,y \in \rr {n} ,
\end{equation}
for some constant $C>0$, and if $v$ in \eqref{moderate} can be chosen
as a polynomial, then $\omega$ is called polynomially
moderated. Furthermore, $v$ is called \emph{submultiplicative} if
\eqref{moderate} holds for $\omega =v$. We denote by
$\mathscr P(\rr {n})$ the set of all polynomially moderated
functions on $\rr {n}$.

\par

Assume that $p,q\in
[1,\infty ]$, and that $\chi \in \mathscr S(\rr n)\setminus 0$. Then
recall that the \emph{modulation space} $M^{p,q}_{(\omega )}(\rr
n)$ is the set of all $f\in \mathscr S'(\rr n)$ such that
\eqref{modnorm} holds. We note that the definition of
$M^{p,q}_{(\omega )}(\rr n)$ is independent of the choice of window
$\chi$, and that different choices of $\chi$ give rise to equivalent
norms. (See Proposition \ref{p1.4} below.)
For conveniency we set $M^p _{(\omega )}= M^{p,p}_{(\omega
)}$. Furthermore, if $\omega \equiv 1$ we also set
$M^{p,q}=M^{p,q}_{(\omega )}$.

\par

The following proposition is a consequence of well-known facts
in \cite {Fe4, Gc2}. Here and in what follows, we let $p'$
denote the conjugate exponent of $p$, i.{\,}e. $1/p+1/p'=1$ should be
fulfilled.

\par

\begin{prop}\label{p1.4}
Assume that $p,q,p_j,q_j\in [1,\infty ]$ for $j=1,2$, and $\omega
,\omega _1,\omega _2,v\in \mathscr P(\rr {2n})$ are such that $\omega$
is $v$-moderate and $\omega _2\le C\omega _1$ for some constant
$C>0$. Then the following are true:
\begin{enumerate}
\item[{\rm{(1)}}] if $\chi \in M^1_{(v)}(\rr n)\setminus 0$, then $f\in
M^{p,q}_{(v)}(\rr n)$ if and only if \eqref {modnorm} holds,
i.{\,}e. $M^{p,q}_{(\omega )}(\rr n)$ is independent of the choice of
$\chi$. Moreover, $M^{p,q}_{(\omega )}$ is a Banach space under the
norm in \eqref{modnorm}, and different choices of $\chi$ give rise to
equivalent norms;

\vrum

\item[{\rm{(2)}}] if  $p_1\le p_2$ and $q_1\le q_2$  then
$$
\mathscr S(\rr n)\hookrightarrow M^{p_1,q_1}_{(\omega _1)}(\rr
n)\hookrightarrow M^{p_2,q_2}_{(\omega _2)}(\rr n)\hookrightarrow
\mathscr S'(\rr n)\text ;
$$

\vrum

\item[{\rm{(3)}}] the $L^2$ product $( \cdo ,\cdo )$ on $\mathscr
S$ extends to a continuous map from $M^{p,q}_{(\omega )}(\rr
n)\times M^{p'\! ,q'}_{(1/\omega )}(\rr n)$ to $\mathbf C$. On the
other hand, if $\nmm a = \sup \abp {(a,b)}$, where the supremum is
taken over all $b\in \mathscr {S}(\rr n)$ such that
$\nm b{M^{p',q'}_{(1/\omega )}}\le 1$, then $\nmm {\cdot}$ and $\nm
\cdot {M^{p,q}_{(\omega )}}$ are equivalent norms;

\vrum

\item[{\rm{(4)}}] if $p,q<\infty$, then $\mathscr S(\rr n)$ is dense in
$M^{p,q}_{(\omega )}(\rr n)$. The dual space of $M^{p,q}_{(\omega
)}(\rr n)$ can be identified
with $M^{p'\! ,q'}_{(1/\omega )}(\rr n)$, through the form $(\cdo  ,\cdo
)_{L^2}$. Moreover, $\mathscr S(\rr n)$ is weakly dense in $M^{\infty
}_{(\omega )}(\rr n)$.
\end{enumerate}
\end{prop}

\par

Proposition \ref{p1.4}{\,}(1) permits us be rather vague concerning
the choice of $\chi \in  M^1_{(v)}\setminus 0$ in
\eqref{modnorm}. For example, if $C>0$ is a constant and $\Omega$ is a
subset of $\mathscr S'$, then $\nm a{M^{p,q}_{(\omega )}}\le C$ for
every $a\in \Omega$, means that the inequality holds for some choice
of $\chi \in  M^1_{(v)}\setminus 0$ and every $a\in
\Omega$. Evidently, for any other choice
of $\chi \in  M^1_{(v)}\setminus 0$, a similar inequality is true
although $C$ may have to be replaced by a larger constant, if necessary.

\par

It is also convenient to let $\mathcal M^{p,q}_{(\omega )}(\rr {n})$
be the completion of $\mathscr S(\rr n)$
under the norm $\nm \cdot {M^{p,q}_{(\omega )}}$. Then $\mathcal
M^{p,q}_{(\omega )}\subseteq M^{p,q}_{(\omega )}$ with equality if and
only if $p<\infty$ and $q<\infty$. It follows that most of the
properties which are valid for $M^{p,q}_{(\omega )}(\rr n)$, also hold
for $\mathcal M^{p,q}_{(\omega )}(\rr n)$.

\par

We also need to use multiplication properties of modulation
spaces. The proof of the following proposition is omitted since the
result can be found in \cite {Fe4, FG1, To2,To7}.

\par

\begin{prop}\label{multprop}
Assume that $p, p_j,q_j\in [1,\infty ]$ and $\omega _j, v\in \mathscr
P(\rr {2n})$ for $j=0,\dots ,N$ satisfy
$$
\frac 1{p_1}+\cdots +\frac 1{p_N}=\frac 1{p_0},\quad \frac
1{q_1}+\cdots +\frac 1{q_N}=N-1+\frac 1{q_0},
$$
and
$$
\omega _0(x,\xi _1+\cdots +\xi _N)\le C\omega _1(x,\xi _1)\cdots
\omega _N(x,\xi _N),\quad x,\xi _1,\dots \xi _N \in \rr n ,
$$
for some constant $C$. Then $(f_1,\dots ,f_N)\mapsto
f_1\cdots f_N$ from $\mathscr S(\rr n)\times \cdots \times \mathscr
S(\rr n)$ to $\mathscr S(\rr n)$ extends uniquely to a continuous map
from $M^{p_1,q_1}_{(\omega _1)}(\rr n)\times \cdots \times
M^{p_N,q_N}_{(\omega _N)}(\rr n)$ to $M^{p_0,q_0}_{(\omega _0)}(\rr
n)$, and
$$
\nm {f_1\cdots f_N}{M^{p_0,q_0}_{(\omega _0)}}\le C\nm
{f_1}{M^{p_1,q_1}_{(\omega _1)}}\cdots \nm
{f_N}{M^{p_N,q_N}_{(\omega _N)}}
$$
for some constant $C$ which is independent of $f_j\in
M^{p_j,q_j}_{(\omega _j)}(\rr n)$ for $i=1,\dots ,N$.

\par

Furthermore, if  $u_0=0$ when $p<\infty$, $v(x,\xi )=v(\xi )\in \mathscr
P(\rr n)$ is submultiplicative, $f\in M^{p ,1}_{(v)}(\rr n)$, and
$\phi ,\psi$ are entire funcitons on $\mathbf C$ with expansions
$$
\phi (z)=\sum _{k=0}^\infty u_kz^k,\quad \psi (z)=\sum _{k=0}^\infty
|u_k|z^k,
$$
then $\phi (f)\in M^{p ,1}_{(v)}(\rr n)$, and
$$
\nm {\phi (f)}{M^{p ,1}_{(v)}}\le C\, \psi (C\nm f{M^{p,1}_{(v)}}),
$$
for some constant $C$ which is independent of $f\in M^{p,1}_{(v)}(\rr
n)$
\end{prop}

\par

\begin{rem}\label{p1.7}
Assume that $p,q,q_1,q_2\in [1,\infty ]$, $\omega _1\in \mathscr P(\rr
n)$ and that $\omega, v \in \mathscr P(\rr {2n})$ are such that
$\omega$ is $v$-moderate. Then the following properties for modulation
spaces hold:
\begin{enumerate}
\item[{\rm{(1)}}]  if $q_1\le\min (p,p')$, $q_2\ge \max (p,p')$ and
$\omega (x,\xi )=\omega _1(x)$,
then $M^{p,q_1}_{(\omega )}\subseteq L^p_{(\omega _0)}\subseteq
M^{p,q_2}_{(\omega )}$. In particular,
$M^2_{(\omega )}=L^2_{(\omega _0)}$;

\vrum

\item[{\rm{(2)}}] if $\omega (x,\xi )=\omega _1(x)$, then
$M^{p,q}_{(\omega )}(\rr n)\hookrightarrow C(\rr n)$ if and only if
$q=1$;

\vrum

\item[{\rm{(3)}}] $M^{1,\infty}$ is a convolution algebra which
contains all measures on $\rr n$ with bounded mass;

\vrum

\item[{\rm{(4)}}] if $x_0\in \rr n$ and $\omega _0(\xi )=\omega
(x_0,\xi )$, then $M^{p,q}_{(\omega )}\cap \mathscr E' =\mathscr
FL^q_{(\omega _0)}\cap \mathscr E'$. Furthermore, if $B$ is a ball
with radius $r$ and center at $x_0$, then
$$
C^{-1}\nm {\widehat f}{L^q_{(\omega _0)}}\le \nm
f{M^{p,q}_{(\omega )}}\le C\nm {\widehat f}{L^q_{(\omega
_0)}},\quad f\in \mathscr E'(B)
$$
for some constant $C$ which only depends on $r$, $n$, $\omega $ and
the chosen window functions;

\vrum

\item[{\rm{(5)}}] if $\omega (x,\xi )=\omega (\xi ,x)$, then
$M^{p}_{(\omega )}$ is invariant under the Fourier transform. A
similar fact holds for partial Fourier transforms;

\vrum

\item[{\rm{(6)}}] for each $x,\xi\in \rr n$ we have
$$
\Vert e^{i\scal\cdo \xi}f(\cdo -x)   \Vert_{M^{p,q}_{(\omega)}}\le C
v(x,\xi) \Vert f \Vert_{M^{p,q}_{(\omega)}},
$$
for some constant $C$;

\vrum

\item[{\rm{(7)}}] if $\tilde {\omega}(x,\xi)=\omega(x,-\xi)$ then
$f\in M^{p,q}_{(\omega)}$ if and only if $\bar f\in
M^{p,q}_{(\tilde\omega)}$.

\end{enumerate}
(See e.{\,}g. \cite {Fe2, Fe3, Fe4, FG1, FG2, FG4, Gc2, To7}.)
\end{rem}

\par

For future references we note that the constant $C_{r,n}$ is
independent of the center of the ball $B$ in (4) in Remark
\ref{p1.7}.

\par

In our investigations we need the following characterization of
modulation spaces.

\par

\begin{prop}\label{modspchar}
Let $\{ x_\alpha \} _{\alpha \in I}$ be a lattice in $\rr n$,
$B_\alpha =x_\alpha +B$ where $B\subseteq \rr n$ is an open ball, and
assume that $f_\alpha \in \mathscr E'(B_\alpha )$ for every $\alpha
\in I$. Also assume that $p,q\in [1,\infty ]$. Then the following is
true:
\begin{enumerate}
\item if 
\begin{equation}\label{fsum}
f = \sum _{\alpha \in I}f_\alpha \quad \text{and}
\quad
F(\xi ) \equiv \Big (\sum _{\alpha \in i}|\widehat f_\alpha (\xi
)\omega (x_\alpha ,\xi )|^p\Big )^{1/p} \in L^q(\rr n),
\end{equation}
then $f\in M^{p,q}_{(\omega )}$, and $f\mapsto \nm F{L^q}$ defines a
norm on $M^{p,q}_{(\omega )}$ which is equivalent to $\nm \cdo
{M^{p,q}_{(\omega )}}$ in \eqref{modnorm};

\vrum

\item if in addition $\cup _\alpha B_\alpha =\rr n$, $\chi \in C^\infty
_0(B)$ satisfies $\sum _\alpha \chi (\cdo -x_\alpha )=1$, $f\in
M^{p,q}_{(\omega )}(\rr n)$, and $f_\alpha =f\, \chi (\cdo -x_\alpha
)$, then $f_\alpha \in \mathscr E'(B_\alpha )$ and \eqref{fsum} is
fulfilled.
\end{enumerate}
\end{prop}

\par

\begin{proof}
(1) Assume that $\chi \in C_0^\infty (\rr n)\setminus 0$ is
fixed. Since there is a bound of overlapping supports of $f_\alpha$,
we obtain
\begin{multline*}
|\mathscr F(f\chi (\cdo -x))(\xi )\omega (x,\xi )| \le \sum  |\mathscr
F(f_\alpha \chi (\cdo -x ))(\xi )\omega (x,\xi )|
\\[1ex]
\le C \Big (\sum |\mathscr F(f_\alpha \chi (\cdo -x ))(\xi
)\omega (x,\xi )|^p\Big )^{1/p}.
\end{multline*}
for some constant $C$. From the support properties of $\chi$, and the
fact that $\omega $ is $v$-moderate for some $v\in \mathscr P(\rr
{2n})$, it follows for some constant $C$ independent of $\alpha$ we
have
$$
|\mathscr F(f_\alpha \chi (\cdo -x))(\xi )\omega (x,\xi )| \le
C|\mathscr F(f_\alpha \chi (\cdo -x))(\xi )\omega (x_\alpha ,\xi )|.
$$
Hence, for some balls $B'$ and $B'_\alpha =x_\alpha +B'$, we get
\begin{multline*}
\Big (\int |\mathscr F(f\chi (\cdo -x))(\xi )\omega (x,\xi )|^p \,
dx\Big )^{1/p}
\\[1ex]
\le C \Big ( \sum _\alpha \int _{B'_\alpha }|\mathscr F(f_\alpha \chi
(\cdo -x))(\xi )\omega (x_\alpha ,\xi )|^p \, dx\Big )^{1/p}
\\[1ex]
\le C \Big ( \sum _\alpha \int _{B'_\alpha }\big (|\widehat f_\alpha
\omega (x_\alpha ,\cdo )|*|\widehat  \chi v(0,\cdo )|(\xi )\big )^p \,
dx\Big )^{1/p}
\\[1ex]
\le C'' \Big ( \sum _\alpha \big (|\widehat f_\alpha \omega (x_\alpha
,\cdo )|*|\widehat  \chi v(0,\cdo )|(\xi )\big )^p \Big )^{1/p}
\le C'' F*|\widehat \chi v(0,\cdo )|(\xi ),
\end{multline*}
for some constants $C'$ and $C''$. Here we have used Minkowski's
inequality in the last inequality. By applying the $L^q$-norm and
using Young's inequality we get
$$
\nm f{M^{p,q}_{(\omega )}}\le C'' \nm {F*|\widehat \chi v(0,\cdo
)|}{L^q}\le C'' \nm F{L^q}\nm {\widehat \chi v(0,\cdo )}{L^1}.
$$
Since we have assumed that $F\in L^q$, it follows that $\nm
f{M^{p,q}_{(\omega )}}$ is finite. This proves (1).

\par

The assertion (2) follows immediately from the general theory of
modulation spaces. (See e.{\,}g. \cite {GH1, Gc2}.) The proof is
complete.
\end{proof}

\par

Next we discuss (complex) interpolation properties for modulation
spaces. Such properties were carefully investigated in
\cite {Fe4} for classical modulation spaces, and thereafter
extended in several directions in \cite{FG2}, where
interpolation properties for coorbit spaces were established. As a
consequence of \cite{FG2} we have the following proposition.
\par

\begin{prop}\label{interpolmod}
Assume that $0<\theta <1$, $p_j,q_j\in [1,\infty ]$ and that $\omega
_j\in \mathscr P(\rr {2n})$ for $j=0,1,2$ satisfy
$$
\frac 1{p_0} = \frac {1-\theta}{p_1}+\frac {\theta}{p_2},\quad
\frac 1{q_0} = \frac {1-\theta}{q_1}+\frac {\theta}{q_2}\quad
\text{and}\quad \omega _0=\omega _1^{1-\theta }\omega _2^\theta .
$$
Then
$$
(\mathcal M^{p_1,q_1}_{(\omega _1)}(\rr n),\mathcal
M^{p_2,q_2}_{(\omega _2)}(\rr n))_{[\theta ]} = \mathcal
M^{p_0,q_0}_{(\omega _0)}(\rr n).
$$
\end{prop}

\medspace

Next we recall some facts in Chapter XVIII in \cite {H} concerning
pseudo-differential operators. Assume that $a\in \mathscr 
S(\rr {2n})$, and that $t\in \mathbf R$ is fixed. Then the
pseudo-differential operator $a_t(x,D)$ in \eqref{e0.5} is a linear
and continuous operator on $\mathscr S(\rr n)$, as remarked in the
introduction. For general $a\in \mathscr S'(\rr {2n})$, the
pseudo-differential operator $a_t(x,D)$ is defined as the continuous
operator from $\mathscr S(\rr n)$ to $\mathscr S'(\rr n)$ with
distribution kernel
\begin{equation}\label{weylkernel}
K_{t,a}(x,y)=(2\pi )^{-n/2}(\mathscr F_2^{-1}a)((1-t)x+ty,y-x),
\end{equation}
Here $\mathscr F_2F$ is the partial
Fourier transform of $F(x,y)\in \mathscr S'(\rr{2n})$ with respect to
the $y$-variable. This definition makes sense, since
the mappings $\mathscr F_2$ and $F(x,y)\mapsto F((1-t)x+ty,y-x)$ are
homeomorphisms on $\mathscr S'(\rr {2n})$. We also note that this
definition of $a_t(x,D)$ agrees with the operator in \eqref{e0.5} when
$a\in \mathscr S(\rr {2m})$.

\par

Furthermore, for any $t\in \mathbf R$ fixed, it follows from the kernel
theorem by Schwartz that the map $a\mapsto
a_t(x,D)$ is bijective from $\mathscr S'(\rr {2n})$ to $\mathscr
L(\mathscr S(\rr n), \mathscr S'(\rr n))$ (see e.{\,}g. \cite {H}).

\par

In particular, if $a\in \mathscr S'(\rr {2m})$ and $s,t\in
\mathbf R$, then there is a unique $b\in \mathscr S'(\rr {2m})$ such
that $a_s(x,D)=b_t(x,D)$. By straight-forward applications of
Fourier's inversion  formula, it follows that
\begin{equation}\label{pseudorelation}
a_s(x,D)=b_t(x,D) \quad \Leftrightarrow \quad b(x,\xi )=e^{i(t-s)\scal
{D_x}{D_\xi}}a(x,\xi ).
\end{equation}
(Cf. Section 18.5 in \cite{H}.)

\par

We end this section by recalling some facts on Schatten-von Neumann
operators and pseudo-differential operators (cf. the introduction).

\par

For each pairs of Hilbert spaces $\mathscr H_1$ and $\mathscr H_2$,
the set $\mathscr I_p(\mathscr H_1,\mathscr H_2)$ is a Banach space
which increases with $p\in [1,\infty ]$, and if $p<\infty$, then
$\mathscr I_p(\mathscr H_1,\mathscr H_2)$ is contained in
the set of compact operators. Furthermore, $\mathscr I_1(\mathscr
H_1,\mathscr H_2)$,
$\mathscr I_2(\mathscr H_1,\mathscr H_2)$ and $\mathscr I_\infty
(\mathscr H_1,\mathscr H_2)$ agree with the set of
trace-class operators, Hilbert-Schmidt operators and continuous
operators respectively, with the same norms.

\par

Next we discuss complex interpolation properties of Schatten-von
Neumann classes. Let $p,p_1,p_2\in [1,\infty ]$ and let $0\le \theta
\le 1$. Then it holds
\begin{equation}\label{interpschatt}
\mathscr I_p =(\mathscr I_{p_1},\mathscr I_{p_2})_{[\theta ]},\quad
\text{when}\quad \frac 1p = \frac {1-\theta}{p_1}+\frac \theta {p_2}.
\end{equation}
We refer to
\cite{Si, To8} for a brief discussion of Schatten-von Neumann
operators.

\par

For any $t\in \mathbf R$ and $p\in [1,\infty ]$, let $s_{t,p}(\omega
_1,\omega _2)$ be the set of all $a\in \mathscr S'(\rr {2n})$ such
that $a_t(x,D)\in \mathscr I_p(M^2_{(\omega _1)},M^2_{(\omega
_2)})$. Also set
$$
\nm a{s_{t,p}}=\nm a{s_{t,p}(\omega _1,\omega _2)}\equiv \nm
{a_t(x,D)}{\mathscr I_p(M^2_{(\omega _1)},M^2_{(\omega _2)})}
$$
when $a_t(x,D)$ is continuous from $M^2_{(\omega _1)}$ to
$M^2_{(\omega _2)}$. By using the fact that $a\mapsto
a_t(x,D)$ is a bijective map from $\mathscr S'(\rr {2n})$ to $\mathscr
L(\mathscr S(\rr n),\mathscr S'(\rr n))$, it follows that the map
$a\mapsto a_t(x,D)$ restricts to an isometric bijection from
$s_{t,p}(\omega _1,\omega _2)$ to $\mathscr I_p(M^2_{(\omega
_1)},M^2_{(\omega _2)})$.

\par

Here and in what follows we let $p'\in [1,\infty ]$ denote the conjugate
exponent of $p\in [1,\infty ]$, i.{\,}e. $1/p+1/p'=1$.

\par

\begin{prop}\label{pseudomod}
Assume that $p,q_1,q_2\in [1,\infty ]$ are such that $q_1\le \min
(p,p')$ and $q_2\ge \max (p,p')$. Also assume that $\omega _1,\omega
_2\in \mathscr P(\rr {2n})$ and $\omega ,\omega _0\in \mathscr P(\rr
{4n})$ satisfy
$$
\frac{\omega _2(x-ty,\xi + (1-t)\eta )}{\omega _1(x+(1-t)y,\xi -t\eta
)}=\omega (x,\xi ,\eta ,y)
$$
and
$$
\omega_0(x,y,\xi,\eta)=\omega((1-t)x+ty, t\xi-(1-t)\eta,\xi+\eta,y-x).
$$
Then the following is true:
\begin{enumerate}
\item $M^{p,q_1}_{(\omega)}(\rr {2n})\subseteq
s_{t,p}(\omega_1,\omega_2)\subseteq M^{p,q_2}_{(\omega)}(\rr {2n})$;

\vrum

\item the operator kernel $K$ of $a_t(x,D)$ belongs to
$M^p_{(\omega_0)}(\rr {2n})$ if
and only if $a\in M^p_{(\omega)}(\rr {2n})$ and for some constant $C$,
which only depends on $t$ and the involved weight functions, it holds
$\Vert K \Vert_{M^p_{(\omega_0)}}=C \Vert a \Vert_{M^p_{(\omega)}}$
\end{enumerate}
\end{prop}

\par

\begin{proof}
The assertion (1) is a restatement of Theorem 4.13 in \cite{Toft4}. The
assertion (2) follows by similar arguments as in the proof of
Proposition 4.8 in \cite{Toft4}, which we recall here. Let  $\chi ,
\psi\in \mathscr S(\rr {2n})$ be such that
$$
\psi(x,y) = \int \chi ((1-t)x+ty,\xi)e^{i\scal {y-x}\xi}\, d\xi.
$$
By applying the Fourier inversion formula it follows by
straightforward computations that
$$
|\mathscr F(K \tau _{(x-ty, x+(1-t)y)}\psi )(\xi +(1-t)\eta ,
-\xi +t\eta )|=|\mathscr F(a\tau _{(x,\xi)}\chi )(y,\eta )|.
$$
The result now follows by applying the $L^p_{(\omega)}$ norm on these
expressions.
\end{proof}

\par

\section{Continuity properties of Fourier
integral operators}\label{sec2}

\par

In this section we extend Theorem 3.2 in \cite{Bu0} in such way that
more general modulation spaces are involved. In these investigations
we assume that  the phase function $\fy$ and the amplitude $a$ depend
on $x, y\in \rr n$ and $\zeta \in \rr m$. For conveniency we use the
notation $X, Y, Z,\dots$ for tripples of the form $(x,y,\zeta )\in \rr
{2n+m}$. We start to make an appropriate definition of the involved
Fourier integral operators.

\par

Assume that $v\in \mathscr P(\rr {2n+m}\times \rr {2n+m})$ is
sub-multiplicative and satisfies
\begin{equation}\label{vvrel}
\begin{aligned}
v(X,\xi ,\eta ,z)&=v(\xi ,\eta ,z),\qquad \text{and}
\\[1ex]
v(t\cdo )\le Cv \qquad \xi ,\eta \in \rr n\ z\in \rr m,
\end{aligned}
\end{equation}
for some constant $C$ which is independent of $t\in [0,1]$
(i.{\,}e. $v(X,\xi ,\eta ,z)$ is constant with respect to $X\in \rr
{2n+m}$. For each real-valued $\fy \in C^{2}(\rr {2n+m})$ which
satisfies $\partial ^\alpha \fy \in M^{\infty ,1}_{(v)}$, and $a\in
\mathscr S(\rr {2n+m})$, it follows that the Fourier integral operator
$f\mapsto \op _\fy (a)f$ in \eqref{fintop} is well-defined and makes
sense as a continuous operator from $\mathscr S(\rr n)$ to $\mathscr
S'(\rr n)$. If $f,g\in \mathscr S(\rr n)$, then
$$
(\op _\fy (a)f,g)=(2\pi )^{-n}\int a(X)e^{i\fy
(X)}f(y)\overline{g(x)}\, dX.
$$
In order to extend the definition we reformulate the latter relation
in terms of short-time Fourier transforms.

\par

Assume that $0\le \chi ,\psi \in
C_0^\infty (\rr {2n+m})$ and $0\le \chi _0\in C_0^\infty (\rr n)$ are
such that
$$
\nm {\chi _0} {L^1}=\nm {\chi }{L^2}=1,
$$
and that $X_1=(x_1,y_1,\zeta _1)\in \rr {2n+m}$ as usual. By
straight-forward computations we get
\begin{multline*}
(\op _\fy (a)f,g) = \int a(X)f(y)\overline {g(x)} e^{i\fy (X )}\, dX
\\[1ex]
= \iint a(X+X_1)\chi (X_1)^2f(y+y_1)\chi _0(y_1)\overline
{g(x+x_1)\chi _0(x_1)} e^{i\psi (X_1)\fy (X+X_1 )}\, dXdX_1
\end{multline*}
If $\mathscr F _{1,2}a$ denotes the partial Fourier transform of
$a(x,y,\zeta )$ with respect to the $x$ and $y$ variables, then
Parseval's formula gives
\begin{multline*}
(\opp _\fy (a)f,g)
\\[1ex]
=\iiiint  F(X,\xi ,\eta ,\zeta _1) \mathscr F(f(y+\cdo )\chi _0)(-\eta
)\overline {\mathscr F(g(x+\cdo )\chi _0)(\xi )}\, dXd\xi d\eta d\zeta
_1
\\[1ex]
=\iiiint  F(X,\xi ,\eta ,\zeta _1) (V_{\chi _0}f)(y,-\eta )\overline
{(V_{\chi _0}g)(x,\xi )}e^{-i(\scal x\xi+\scal y{\eta})}\, dXd\xi
d\eta d\zeta _1,
\\[1ex]
=\iiint \Big ( \int  F(X,\xi ,\eta ,\zeta _1)\, d\zeta _1 \Big )
(V_{\chi _0}f)(y,-\eta )\overline {(V_{\chi _0}g)(x,\xi )}e^{-i(\scal
x\xi+\scal y{\eta})}\, dXd\xi d\eta ,
\end{multline*}
where
\begin{equation*}
F(X,\xi ,\eta ,\zeta _1) = \mathscr F_{1,2}\big (e^{i\psi (\cdo ,\zeta
_1)\fy (X+(\cdo ,\zeta _1 ))}a(X+(\cdo ,\zeta _1))\chi (\cdo ,\zeta
_1)^2\big )(\xi ,\eta ).
\end{equation*}

\par

By Taylor's formula it follows that
$$
\psi (X_1) \fy (X+X_1) =\psi (X_1)\psi _{1,X}(X_1) +\psi _2{2,X}(X_1),
$$
where
\begin{equation}\label{psidef}
\begin{aligned}
\psi _{1, X}( X_1) &= \fy ( X)+\scal {\fy
'( X)}{ X_1}
\\[1ex]
\psi _{2, X}( X_1) &= \psi ( X_1)\int _0^1
(1-t)\langle \varphi ''( X+t X_1)
X_1, X_1\rangle \, dt .
\end{aligned}
\end{equation}
By inserting these expressions into the definition of $F(X,\xi ,\eta
,\zeta _1)$, and integrating with respect to the $\zeta _1$-variable
give
\begin{multline*}
\int F(X,\xi ,\eta ,\zeta _1)\, d\zeta _1
\\[1ex]
=\mathscr F((e^{i\psi _{2,X}}\chi ) (a(\cdo +X)\chi )(\xi -\fy
'_x(X),\eta -\fy '_y(X),-\fy '_\zeta (X))
\\[1ex]
=\mathcal H_{a,\fy}(X,\xi ,\eta ),
\end{multline*}
where
\begin{equation}\label{Hdef}
\begin{aligned}
\mathcal H_{a,\fy}(X,\xi ,\eta )
= h_X*(\mathscr F(a(\cdo +&X)\chi ))(\xi -\fy '_x(X),\eta -\fy
'_y(X),-\fy '_\zeta (X)),
\\[1ex]
\text{and}\quad h_X &= (2\pi )^{-n} (\mathscr F(e^{i\psi _{2,X}}\chi ))
\end{aligned}
\end{equation}

\par

Summing up we have proved that
\begin{equation}\label{fourrel}
\begin{aligned}
(&\operatorname{Op} _{\fy} (a)f,g) =T_{a,\fy}(f,g)
\\[1ex]
&\equiv \iiint  \mathcal
H_{a,\fy}(X,\xi ,\eta )
(V_{\chi _0}f)(y,-\eta )\overline {(V_{\chi _0}g)(x,\xi )}e^{-i(\scal
x\xi+\scal y{\eta})}\, dXd\xi d\eta .
\end{aligned}
\end{equation}

\par

If $a\in \mathscr S'(\rr {2n+m})$, $f,g\in \mathscr S(\rr n)$ and that
the mapping
$$
(X,\xi ,\eta )\mapsto \mathcal H_{a,\fy}(X,\xi ,\eta ) (V_{\chi
_0}f)(y,-\eta )\overline {(V_{\chi _0}g)(x,\xi )}
$$
belongs to $L^1(\rr {2n+m}\times \rr {2n})$, then we still let
$T_{a,\fy}(f,g)$ be defined as the right-hand side of
\eqref{fourrel}. In what follows we use \eqref{fourrel} to extend the
definition of Fourier integral operator with more general
amplitudes. Here recall that if for each fixed $f_0\in \mathscr S$ and
$g_0\in \mathscr S$, the mappings $f\mapsto T(f,g_0)$ and $g\mapsto
T(f_0,g)$ are continuous from $\mathscr S$ to $\mathbf C$, then it
follows by Banach-Steinhauss theorem that
$$
(f,g)\mapsto T(f,g)
$$
is continuous from $\mathscr S\times \mathscr S$ to $\mathbf C$.

\par

\begin{defn}\label{deffourop}
Assume that $v\in \mathscr P(\rr {2n+m}\times \rr {2n+m})$ is
submultiplicative and satisfies \eqref{vvrel}, $\fy \in C^2(\rr
{2n+m})$ is such and that $\partial ^\alpha \fy \in M^{\infty
,1}_{(v)}$, and that $a\in \mathscr S'(\rr {2n+m})$ is such that
$f\mapsto T_{a,\fy}(f,g_0)$ and $g\mapsto T_{a,\fy}(f_0,g)$ are
well-defined and continuous from $\mathscr S(\rr n)$ to $\mathbf C$,
for each fixed $f_0,g_0\in \mathscr S(\rr n)$. Then $(a,\fy )$ is
called \emph{admissible}, and the Fourier integral operator $\op _\fy
(a)$ is the continuous mapping from $\mathscr S(\rr n)$ to $\mathscr
S'(\rr n)$ which is defined by the formulas \eqref{psidef},
\eqref{Hdef} and \eqref{fourrel}.
\end{defn}

\par

Our general continuity results also includes weights which satisfy
conditions of the form
\begin{equation}\label{weightsineq}
\begin{aligned}
\frac {\omega _2(x,\xi )}{\omega _1(y,-\eta )} \le C\omega (X,\xi -\fy
'_x(X) &, \eta -\fy '_y(X),-\fy '_\zeta (X)),
\\[1ex]
X &= (x,y,\zeta )\in \rr {2n+m},\quad \xi ,\eta \in \rr n.
\end{aligned}
\end{equation}

\par

\begin{thm}\label{thm3.1A}
Assume that $\mathsf d>0$, $\omega ,v\in \mathscr P(\rr {2n+m}\times
\rr {2n+m})$ and $\omega _1,\omega _2\in \mathscr P(\rr {2n})$ are
such that \eqref{vvrel} and \eqref{weightsineq} are fulfilled and that
$v$ is submultiplicative. Also assume that $\fy \in C^2(\rr {2n+m})$
is such that $\partial ^\alpha \fy \in M^{\infty ,1}_{(v)}$,
$a\in \mathscr S'(\rr {2n+m})$, and that one of the following
conditions hold:
\begin{enumerate}
\item $|\det (\fy ''_{\zeta ,\zeta })|\ge \mathsf d$ and $\nmm
a<\infty$, where
\begin{equation}\label{thm3.1Aineq1}
\nmm a = \sup _{X,\xi ,\eta}\Big (\int |V_\chi a(X,\xi ,\eta ,z)\omega
(X,\xi ,\eta ,z)|\, dz\Big ) <\infty ;
\end{equation}

\vrum

\item $m=n$, $|\det (\fy ''_{x,\zeta })|\ge \mathsf d$ and $\nmm
a<\infty$, where
\begin{equation}\label{thm3.1Aineq2}
\nmm a = \sup _{X,\eta ,z}\Big (\int |V_\chi a(X,\xi ,\eta ,z)\omega
(X,\xi ,\eta ,z)|\, d\xi\Big ) \text ;
\end{equation}

\vrum

\item $m=n$, $|\det (\fy ''_{y,\zeta })|\ge \mathsf d$ and $\nmm
a<\infty$, where
\begin{equation}\label{thm3.1Aineq3}
\nmm a = \sup _{X,\xi ,z}\Big (\int |V_\chi a(X,\xi ,\eta ,z)\omega
(X,\xi ,\eta ,z)|\, d\eta \Big ) .
\end{equation}
\end{enumerate}
Then $(a,\fy )$ is admissible, and $\op _\fy (a)$ extends to a linear
and continuous operator from $M^1_{(\omega _1)}(\rr n)$ to $M^\infty
_{(\omega _2)}(\rr n)$.

\par

Moreover, for some $C$ which is independent of $a$ and $\fy$ it holds
\begin{equation}\label{opest1}
\nm {\op _\fy (a)}{M^1_{(\omega _1)}\to M^\infty _{(\omega _2)}} \le
\frac {C\nmm a}{\mathsf d}\exp (C\nm {\fy ''}{M^{\infty ,1}_{(v)}}).
\end{equation}
\end{thm}

\par

We need some preparing lemmas for the proof.

\par

\begin{lemma}\label{lemmafourop3}
Assume that $v(x,\xi )=v(\xi )\in \mathscr P(\rr {n})$ is
submultiplicative and satisfies $v(t\xi )\le Cv(\xi )$ for some
constant $C$ which is independent of $t\in [0,1]$ and $\xi \in \rr
n$. Also assume that $f\in M^{\infty ,1}_{(v)}(\rr n)$, $\chi \in
C^\infty _0(\rr n)$ and that $x\in \rr n$, and let
$$
\fy _{x,j,k}(y)=\chi (y)\int _0^1 (1-t)f(x+ty)y_jy_k\, dt .
$$
Then there is a constant $C$ and a function $g\in M^1_{(v)}(\rr n)$
such that $\nm {g}{M^1_{(v)}}\le C\nm f{M^{\infty ,1}_{(v)}}$ and
$|\mathscr F(\fy _{x,j,k})(\xi )|\le \widehat g(\xi )$.
\end{lemma}

\par

\begin{proof}
We first prove the assertion when $\chi$ is replaced by $\psi
_0(y)=e^{-2|y|^2}$. For conveniency we let 
$$
\psi _1(y)= e^{-|y|^2},\quad \text{and}\quad \psi _2(y) = e^{|y|^2},
$$
and 
$$
H_{\infty ,f} (\xi )\equiv \sup _{x}|\mathscr F(f\, \psi _1(\cdo
-x))(\xi )|.
$$
We claim that $g$, defined by
\begin{equation}\label{defg}
\widehat g(\xi)=\int_0^1\int (1-t)H_{\infty,f}(\eta)e^{-|\xi-
t\eta|^2/16}\,d\eta dt,
\end{equation}
fulfills the required properties.

\par

In fact, by applying $M^1_{(v)}$ norm on $g$ and using Minkowski's
inequality and Remark \ref{p1.7} (6), we obtain
\begin{align*}
\nm g{M^1_{(v)}} &= \Big \Vert \int_0^1\int
(1-t)H_{\infty,f}(\eta)e^{-|\xi- t\eta|^2/16}\,d\eta dt \Big \Vert
_{M^1_{(v)}}
\\[1ex]
&\le \int_0^1\int (1-t)H_{\infty,f}(\eta)\nm {e^{-|\cdo -
t\eta|^2/16}}{M^1_{(v)}}\,d\eta dt
\\[1ex]
&\le C_1\int_0^1\int (1-t)H_{\infty,f}(\eta)\nm {e^{-|\cdo
|^2/16}}{M^1_{(v)}}v(t\eta )\,d\eta dt
\\[1ex]
&\le C_2\int_0^1\int (1-t)H_{\infty,f}(\eta)v(\eta )\nm {e^{-|\cdo
|^2/16}}{M^1_{(v)}}\,d\eta dt
\\[1ex]
&= C_3\nm {H_{\infty,f}v}{L^1} = C_3\nm f{M^{\infty ,1}_{(v)}}.
\end{align*}

\par

In order to prove that $|\mathscr F(\fy _{x,j,k})(\xi )|\le g(\xi )$,
we let $\psi (y)=\psi _{j,k}(y)= y_jy_k\psi _0 (y)$. Then
$$
\fy _{x,j,k}(y)= \psi (y)\int _0^1 (1-t)f(x+ty)\, dt.
$$
By a change of variables we obtain
\begin{equation}\label{estfop3}
\begin{aligned}
|\mathscr F(\fy _{x,j,k})(\xi )| &= \Big |\int _0^1 (1-t)\Big ( \int
f(x+ty)\psi (y)e^{-i\scal y\xi}\, dy\Big )\, dt \Big |
\\[1ex]
&= \Big |\int _0^1 t^{-n}(1-t)\mathscr F(f\, \psi ((\cdo -x)/t))(\xi
/t)e^{i\scal x\xi /t}\, dt\Big |
\\[1ex]
&\le \int _0^1 t^{-n}(1-t)\sup _{x\in \rr n}|\mathscr F(f\,
\psi ((\cdo -x)/t))(\xi /t)|\, dt.
\end{aligned}
\end{equation}
We need to estimate the right-hand side. By straight-forward
computations we get
\begin{multline*}
|\mathscr F(f\, \psi ((\cdo -x)/t))(\xi )|
\\[1ex]
\le  (2\pi )^{-n/2}\big
( |\mathscr F(f\, \psi _1(\cdo -x))|*| \mathscr F(\psi
((\cdo -x)/t)\, \psi _2(\cdo -x))|\big )(\xi )
\\[1ex]
=(2\pi )^{-n/2}\big
( |\mathscr F(f\, \psi _1(\cdo -x))|*| \mathscr F(\psi
(\cdo /t)\, \psi _2)|\big )(\xi )
\end{multline*}
where the convolutions should be taken with respect to the
$\xi$-variable only. This implies that
\begin{equation}\label{estfop33}
|\mathscr F(f\, \psi ((\cdo -x)/t))(\xi )|
\le (2\pi )^{-n/2}\big
( H_{\infty ,f}*| \mathscr F(\psi (\cdo /t)\, \psi _2)|\big )(\xi )
\end{equation}

\par

In order to estimate the latter Fourier transform we note that
\begin{equation}\label{gaussest1}
| \mathscr F(\psi (\cdo /t)\, \psi _2)| =  |\partial _j\partial _k
\mathscr F(\psi _0(\cdo /t)\, \psi _2)|.
\end{equation}
Since $\psi _0$ and $\psi _2$ are Gauss functions and $0\le t\le 1$, a
straight-forward computation gives
\begin{equation}\label{gaussest2}
\mathscr F(\psi _0(\cdo /t)\, \psi _2)(\xi ) = \pi
^{n/2}t^n(2-t^2)^{-n/2}e^{-t^2|\xi |^2/(4(2-t^2))}.
\end{equation}
A combination of \eqref{gaussest1} and \eqref{gaussest2} therefore
give
\begin{equation}\label{gaussest3}
| \mathscr F(\psi (\cdo /t)\, \psi _2)(\xi )|\le
Ct^ne^{-t^2|\xi |^2/16},
\end{equation}
for some constant $C$ which is independent of $t\in [0,1]$. The
assertion now follows by combining \eqref{estfop3}, \eqref{estfop33}
and \eqref{gaussest3}.

\par

In order  to prove the result for general $\chi\in C^\infty_0(\rr n)$
we set
$$
h_{x,j,h}(y)=\psi_0(y)\int_0^1 (1-t)f(x+ty)y_jy_k\, dt,
$$
and we observe that the result is already proved  when
$\varphi_{x,j,k}$ is replaced by $h_{x,j,h}$ and moreover $\varphi
_{x,j,k}=\chi_1 h_{x,j,k}$, for some $\chi_1\in C^\infty_0(\rr n)$.
Hence if $g_0$ is given as the right-hand side  of \eqref{defg}, the
first part of the proof  shows that
$$
\left| \mathscr F(\varphi_{x,j,k})(\xi)\right|=\left| \mathscr
F(\chi_1 h_{x,j,k}(\xi))\right|\le (2\pi)^{-n/2}|\widehat{\chi_1}|\ast
g-0(\xi)\equiv g(\xi).
$$
Moreover $\Vert g_0 \Vert_{M^1_{(v)}} \le C\Vert f \Vert_{M^{\infty,
1}_{(v)}}$.

\par

Since $M^1_{(v)}\ast L^1_{(v)}\subseteq M^1_{(v)}$, we get for some
positive constants $C, C_1$
$$
\Vert g \Vert_{M^1_{(v)}}\le C\Vert
\widehat{\chi_1}\Vert_{L^1_{(v)}}\Vert g_0 \Vert_{M^1_{(v)}}\le C_1
\Vert f \Vert_{M^{\infty, 1}_{(v)}},
$$
which proves the result
\end{proof}

\par

As a consequence of Lemma \ref{lemmafourop3} we have the following
result.

\par

\begin{lemma}\label{lemmafourop33}
Assume that $v(x,\xi )=v(\xi )\in \mathscr P(\rr {n})$ is
submultiplicative and satisfies $v(t\xi )\le Cv(\xi )$ for some
constant $C$ which is independent of $t\in [0,1]$ and $\xi \in \rr
n$. Also assume that $f_{j,k}\in M^{\infty ,1}_{(v)}(\rr n)$ for
$j,k=1,\dots, n$, $\chi \in
C^\infty _0(\rr n)$ and that $x\in \rr n$, and let
$$
\fy _{x}(y)=\sum _{j,k=1,\dots ,n}\fy _{x,j,k}(y),\quad
\text{where}\quad \fy _{x,j,k}(y)=\chi (y)\int _0^1
(1-t)f_{j,k}(x+ty)y_jy_k\, dt .
$$
Then there is a constant $C$ and a function $\Psi \in M^1_{(v)}(\rr
n)$ such that
$$
\nm {\Psi}{M^1_{(v)}}\le \exp(C\sup _{j,k}\nm {f_{j,k}}{M^{\infty
,1}_{(v)}})
$$
and
\begin{equation}\label{psiestim}
|\mathscr F(\exp (i\fy _{x})(\xi ))|\le (2\pi )^{n/2}\delta _0
+\widehat \Psi (\xi ).
\end{equation}
\end{lemma}

\par

\begin{proof}
By Lemma \ref{lemmafourop3}, we may find a function $g\in M^1_{(v)}$
and a constant $C$ such that
$$
|\widehat {\fy _x}(\xi )| \le \widehat g(\xi ),\qquad \nm
g{M^1_{(v)}}\le C \sup _{j,k}(\nm {f_{j,k}}{M^{\infty ,1}_{(v)}}.
$$
Set
\begin{equation*}
\begin{alignedat}{3}
\Phi _{0,x} &= (2\pi )^{n/2}\delta _0,&\qquad \Phi _{l,x} &= |\mathscr
F(\fy _{x})| * \cdots *|\mathscr F(\fy _{x})|,&\qquad l&\ge 1
\\[1ex]
\Upsilon _0 &= (2\pi )^{n/2}\delta _0,&\qquad \Upsilon _l &= g*\cdots
*g,&\qquad l&\ge 1,
\end{alignedat}
\end{equation*}
with $l$ factors in the convolutions. Then by Taylor series expanding,
there is a constant $C$ such that
\begin{equation*}
|\mathscr F(\exp (i\fy _{x}(\cdo ))(\xi ))| \le \sum _{l=0}^\infty
C^{l}\Phi _{l,x}/l! \le  \sum _{l=0}^\infty C^{l}\Upsilon _l/l!
\end{equation*}
Hence, if we set
$$
\Psi =  \sum _{l=1}^\infty C^{l}\Upsilon _l/l!,
$$
it follows that \eqref{psiestim} holds.

\par

Furthermore, since $v$ is submultiplicative we have
$$
\nm {\Upsilon _l}{M^{1}_{(v)}} = \nm {g*\cdots *g}{M^{1}_{(v)}} \le
(C_1\nm g{M^{1}_{(v)}})^l,
$$
for some constant $C_1$. This gives
\begin{multline*}
\nm \Psi{M^1_{(v)}}\le \sum _{l=1}^\infty \nm {\Upsilon
_l}{M^{1}_{(v)}}
\\[1ex]
\le \sum _{l=1}^\infty (C_1\nm g{M^{1}_{(v)}})^l \le \sum
_{l=1}^\infty (C_2\sup _{j,k}(\nm {f_{j,k}}{M^{\infty
,1}_{(v)}})^l
\le \exp (C_2\sup _{j,k}\nm {f_{j,k}}{M^{\infty ,1}_{(v)}}),
\end{multline*}
for some constants $C_1$ and $C_2$. This proves the assertion.
\end{proof}

\par

\begin{proof}[Proof of Theorem \ref{thm3.1A}]
We only prove the result when (1) if fulfilled. The other cases follow
by similar arguments and are left for the reader.

\par

Since
$$
|\mathscr F(e^{i\Psi _{2,X}}\chi )|\le (2\pi )^{-n+m/2}|\mathscr
F(e^{i\Psi _{2,X}}|*|\widehat \chi |,\qquad |\mathscr F(a(\cdo +X)\chi
)| = |V_\chi a(X,\cdo )|,
$$
\eqref{Hdef},  and Lemma \ref{lemmafourop33} give
$$
|\mathcal H_{a,\fy}(X,\xi ,\eta )|\le C (G*|V_\chi a(X,\cdo ))(\xi
-\fy '_x(X), \eta -\fy '_y(X),-\fy '_\zeta (X))|,
$$
for some $G\in L^1_(v)$ which satisfies $\nm G{L^1_{(v)}}\le C\exp
(C\nm {\fy ''}{M^{\infty ,1}_{(v)}}$. By combining this with
\eqref{weightsineq} and letting
\begin{equation}\label{Efuncdef}
\begin{aligned}
E_{a,\omega}(\xi ,\eta ,z) &= \sup _{X}|V_\chi a(X,\xi ,\eta ,z ))
\omega (X,\xi ,\eta ,z )|
\\[1ex]
F_1(x,\xi ) &= |V_{\chi _0}f(x,\xi )\omega _1(x,\xi )|,
\\[1ex]
F_2(x,\xi ) &= |V_{\chi _0}g(x,\xi )/\omega _2(x,\xi )|,
\end{aligned}
\end{equation}
we get
\begin{multline*}
\iiint | \mathcal H_{a,\fy}(X,\xi ,\eta ) (V_{\chi _0}f)(y,-\eta
){(V_{\chi _0}g)(x,\xi )} |\, dXd\xi d\eta
\\[1ex]
\le
C_1\iiint (G*|V_\chi a(X,\cdo )|)(\xi -\fy '_x(X), \eta -\fy
'_y(X),-\fy '_\zeta (X))\times
\\[1ex]
|(V_{\chi _0}f)(y,-\eta ){(V_{\chi _0}g)(x,\xi )}|\, dXd\xi d\eta
\\[1ex]
\le
C_2\iiint (G*E_{a,\omega }(\xi -\fy
'_x(X) , \eta -\fy '_y(X),-\fy '_\zeta (X))\times
\\[1ex]
F_1(y,-\eta )F_2(x,\xi )\, dXd\xi d\eta
\\[1ex]
\le
C_3\nm G{L^1_{(v)}}\iiint (E_{a,\omega} (\xi -\fy '_x(X) , \eta -\fy
'_y(X),-\fy '_\zeta (X))\times
\\[1ex]
F_1(y,-\eta )F_2(x,\xi )\, dXd\xi d\eta .
\end{multline*}
Summing up we have proved that
\begin{equation}\label{estfourop1}
\begin{aligned}
&|(\op _\fy (a)f,g)|\le C\exp (C\nm {\fy ''}{M^{\infty ,1}_{(v)}})\times
\\[1ex]
&\times \iiint (E_{a,\omega} (\xi -\fy '_x(X) , \eta -\fy
'_y(X),-\fy '_\zeta (X)) F_1(y,-\eta )F_2(x,\xi )\, dXd\xi d\eta .
\end{aligned}
\end{equation}
By taking $x,y,-\fy '_\zeta (X),\xi ,\eta $ as new variables of
integration, and using the fact that $|\det (\fy ''_{\zeta
,\zeta})|\ge \mathsf d$ we get
\begin{multline*}
\iiint | \mathcal H_{a,\fy}(X,\xi ,\eta ) (V_{\chi _0}f)(y,-\eta
){(V_{\chi _0}g)(x,\xi )} |\, dXd\xi d\eta
\\[1ex]
\le
\frac {C_3}{\mathsf d} \nm {\fy ''}{M^{\infty ,1}_{(v)}} \iiiint \Big
(  \int E_{a,\omega}(\xi -\fy '_x(X) , \eta -\fy '_y(X),z)\, dz\Big
)\times
\\[1ex]
|(V_{\chi _0}f)(y,-\eta )\omega _1(y,-\eta ){(V_{\chi _0}g)(x,\xi
)}/\omega _2(x,\xi )|\, dxdyd\xi d\eta
\\[1ex]
\le
\frac {C_3\nmm a}{\mathsf d} \exp (C\nm {\fy ''}{M^{\infty ,1}_{(v)}})
\iiiint |(V_{\chi _0}f)(y,-\eta )\omega _1(y,-\eta ){(V_{\chi
_0}g)(x,\xi )}/\omega _2(x,\xi )|\, dxdyd\xi d\eta
\\[1ex]
= \frac {C_3\nmm a}{\mathsf d} \exp (C\nm {\fy ''}{M^{\infty
,1}_{(v)}}) \nm f{M^1_{(\omega _1)}}\nm g{M^1_{(1/\omega _2)}}.
\end{multline*}
This proves that \eqref{opest1} holds, and the result follows.
\end{proof}

\par

Next we consider Fourier integral operators with symbols in $M^{\infty
,1}_{(\omega )}(\rr {2n+m})$. The following result generalizes Theorem
3.2 in \cite {Bu1}.

\par

\begin{thm}\label{boulkemA}
Assume that $1<p<\infty$, $\mathsf d >0$, $\omega ,v\in \mathscr P(\rr
{2n+m}\times \rr {2n+m})$ and $\omega _1,\omega _2\in \mathscr P(\rr
{2n})$ are such that \eqref{vvrel} and \eqref{weightsineq} are
fulfilled and that $v$ is submultiplicative. Also assume that $\fy \in
C^2(\rr {2n+m})$ is such that $\partial ^\alpha \fy \in M^{\infty
,1}_{(v)}$ for $|\alpha | =2$ and \eqref{detphicond} are
fulfilled. Then the following is true:
\begin{enumerate}
\item the map $a\mapsto \op _\fy (a)$ from $\mathscr S(\rr {2n+m})$ to
$\mathscr L(\mathscr S(\rr n),\mathscr S'(\rr n))$ extends uniquely to
a continuous map from $M^{\infty ,1}_{(\omega )}(\rr {2n+m})$ to
$\mathscr L(\mathscr S(\rr n),\mathscr S'(\rr n))$;

\vrum

\item if $a\in M^{\infty ,1}_{(\omega )}(\rr {2n+m})$, then the map
$\op _{\fy }(a)$ from $\mathscr S(\rr n)$ to $\mathscr S'(\rr n)$ is
uniquely extendable to a continuous operator from $M^p_{(\omega
_1)}(\rr n)$ to $M^p_{(\omega _2)}(\rr n)$. Moreover, for some
constant $C$ it holds
\begin{equation}\label{normuppsk2}
\nm {\op _{\fy }(a)}{M^p_{(\omega _1)}\to M^p_{(\omega _2)}}\le
C\mathsf d ^{-1}\nm a{M^{\infty ,1}_{(\omega )}}\exp (C\nm {\fy
''}{M^{\infty 1,}_{(v)}}).
\end{equation}
\end{enumerate}
\end{thm}

\par

\begin{proof}
We shall mainly follow the proof of Theorem 3.2 in \cite{Bu1}. First
assume that $a\in C_0^\infty (\rr {2n+m})$. Let $f,g\in \mathscr S(\rr
n)$. Then it follows tha $\op
_\fy (a)$ makes sense as a continuous operator from $\mathscr S$ to
$\mathscr S'$. By letting
$$
C_\fy =C\exp(C\nm {\fy ''}{M^{\infty ,1}_{(v)}}),
$$
where $C$ is the same as in \eqref{estfourop1}, it follows from
\eqref{Efuncdef}, \eqref{estfourop1} and H{\"o}lder's inequality that
\begin{equation}\label{eq.11}
|(\op _\fy (a)f,g)|\le C_\fy J_1\cdot J_2,
\end{equation}
where
\begin{align*}
J_1 &= \Big ( \iiint (E_{a,\omega} (\xi -\fy '_x(X) , \eta -\fy
'_y(X),-\fy '_\zeta (X)) F_1(y,-\eta )^p\, dXd\xi d\eta \Big )^{1/p}
\\[1ex]
J_2 &= \Big ( \iiint (E_{a,\omega} (\xi -\fy '_x(X) , \eta -\fy
'_y(X),-\fy '_\zeta (X)) F_2(y,-\eta )^{p'}\, dXd\xi d\eta \Big
)^{1/p'}.
\end{align*}
We have to estimate $J_1$ and $J_2$. By taking $z =\fy '_3(X)$,
$\zeta _0=\fy '_2(X)$, $y$, $\xi$ and $\eta$ as new variables of
integrations, and using \eqref{detphicond}, it follows that
\begin{multline*}
J_1 \le  \Big (  \mathsf d ^{-1}\iiint (E_{a,\omega} (\xi -\fy '_x(X)
, \eta -\zeta _0,z) F_1(y,-\eta )^p\, dy dz d\xi d\eta d\zeta _0\Big
)^{1/p}
\\[1ex]
= \Big (  \mathsf d ^{-1}\iiint (E_{a,\omega} (\xi, \zeta _0,z)
F_1(y,-\eta )^p\, dy dz d\xi d\eta d\zeta _0\Big )^{1/p}
\\[1ex]
= \mathsf d ^{-1/p}\nm {E_{a,\omega }}{L^1}^{1/p}\nm {F_1}{L^p}.
\end{multline*}
Hence
\begin{equation}\label{eq.22}
J_1\le \mathsf d ^{-1/p}\nm {a}{M^{\infty ,1}_{(\omega )}}^{1/p}\nm
{f}{M^p_{(\omega _1)}}.
\end{equation}
If we instead take $x$, $y_0=\fy '_3(X)$, $\xi $, $\eta $ and $\zeta
_0=\fy '_1(X)$ as new variables of integrations, it follows by similar
arguments that
\begin{equation}\tag*{(\ref{eq.22})$'$}
J_2\le \le \mathsf d ^{-1/p'}\nm {a}{M^{\infty ,1}_{(\omega )}}^{1/p'}\nm
{g}{M^{p'}_{(1/\omega _2)}}.
\end{equation}
A combination of \eqref{eq11}, \eqref{eq22} and \eqref{eq22}$'$ now give
$$
|(\op _{\fy}(a)f,g)| \le  C\mathsf d ^{-1}\nm a{M^{\infty
,1}_{(\omega )}} \nm {f}{M^p_{(\omega _1)}} \nm {g}{M^{p'}_{(1/\omega
_2)}}\exp (\nm {\fy ''}{M^{\infty ,1}_{(v)}}),
$$
which proves \eqref{normuppsk2}, and the result follows in this case.

\par

Since $\mathscr S$ is dense in $M^p_{(\omega _1)}$ and
$M^{p'}_{(1/\omega _2)}$, the result also holds for $a\in C_0^\infty $
and $f\in M^{p_{(\omega _1)}}$. Hence it follows by Hahn-Banach's
theorem that the asserted extension of the map $a\mapsto \op _\fy (a)$
exists.

\par

It remains to prove that this extension is unique. Therefore assume
that $a\in M^{\infty ,1}_{(\omega )}$ is arbitrary, and take a
sequence $a_j\in C_0^\infty$ for $j=1,2,\dots$ which converges
to $a$ with respect to the narrow convergence (cf. \cite{Sj1,
Toft4}). Then $E_{a_j,\omega}$ converges to $E_{a,\omega}$ in $L^1$ as
$j$ turns to infinity. By \eqref{psidef}--\eqref{fourrel} and the
arguments at the above, it follows from Lebesgue's theorem that
$$
(\op _\fy (a_j)f,g)\to (\op _\fy (a)f,g)
$$
as $j$ turns to infinity. This proves the uniqueness, and the result
follows.
\end{proof}

\par

\begin{rem}
Assume that $a\in M^{\infty ,1}(\rr {2n+m})$ and that the assumptions
on $\fy$ in Definition \ref{deffourop} is fulfilled with $v\equiv
1$. Also assume that $\kappa \in C_0^\infty (\rr m)$ satisfies $\kappa
(0)=1$. Then it is proved in  \cite{Bu0}, the Fourier integral
operator
\end{rem}

\par

\section{Schatten-von Neumann properties of Fourier
integral operators}\label{sec3}

\par

In this section we discuss Schatten-von Neumann operators for Fourier
integral operators with symbols in $M^{p,q}_{(\omega)}(\rr{2n})$ and
phase functions in $M^{\infty,1}_{(v)}(\rr {3n})$, for appropriate
$\omega$ and $v$.
In these investigations we assume that  the phase functions depend on
$x, y, \xi \in \rr n$ and that the symbols are independent of the $y$
variable, and for conveniency we use the notation $X,
Y, Z,\dots$ for tripples of the form $(x,y,\xi )\in \rr
{3n}$. 

\par

In order to establish a weighted version of Theorem 2.5 in \cite {CT1}
we list some conditions for the weight and phase functions. In what
follows we assume that $\fy \in C^2(\rr {3n})$,  $\omega _0,\omega \in
\mathscr P(\rr {4n})$, $v_1\in \mathscr P(\rr n)$, $v_2\in \mathscr
P(\rr {2n})$, $v\in \mathscr P(\rr {6n})$ and $ s_j,t_j\in
\mathbf R$ for $j=1,2$ satisfy
\begin{equation}\label{phasecond}
|\det(\fy ''_{y,\xi }(x,y,\xi ))|\ge \mathsf d
\end{equation}
\begin{equation}\label{weightcond1}
\omega _0(x,y,\xi ,\varphi '_y(x,y,\eta)) = \omega (x,\eta,\xi-\varphi
'_x(x,y,\eta),-\varphi '_\eta(x,y,\eta))
\end{equation}
and
\begin{equation}\label{weightcond2}
\begin{alignedat}{2}
\omega _0(x,y,\xi , \eta +\zeta ) &\le  \omega _0(x,y,\xi ,\eta
)v_1(\zeta ),&\quad x,y,\xi,\eta ,\zeta &\in \rr n
\\[1ex]
\omega (x,\eta ,\xi _1+\xi _2,y_1+y_2) &\le \omega (x,\eta ,\xi
_1,y_1)v_2(\xi _2,y_2),&\quad x,y_1,y_2, \xi _1,\xi _2,\eta &\in \rr n
\\[1ex]
v(x,y,\zeta,\xi ,\eta ,z) &= v_1(\eta )v_2(\xi ,z),&\quad x,y,z,\xi
,\eta ,\zeta &\in \rr n.
\end{alignedat}
\end{equation}

\par

\begin{thm}\label{fourop3}
Assume that $p\in [1,\infty ]$, $\mathsf d >0$, $v\in \mathscr P(\rr {6n})$ is
submultiplicative, $\omega _0,\omega \in \mathscr P(\rr {4n})$ and
that $\varphi \in C(\mathbf {R}^{3n})$ are such that $\fy$ is
real-valued, $\partial ^{\alpha } \varphi \in M ^{\infty,1}_{(v)}$ for
$|\alpha | =2$ and \eqref{phasecond}--\eqref{weightcond2} are
fulfilled. Then the map
$$
a\mapsto K_{a,\fy}(x,y)\equiv \int a(x, \xi
)e^{i\varphi (x,y,\xi )}\, d\xi ,
$$
from $\mathscr S(\rr {2n})$ to $\mathscr S'(\rr {2n})$ extends
uniquely to a continuous map from $M^p_{(\omega )}(\rr {2n})$ to
$M^p_{(\omega _0)}(\rr {2n})$.
\end{thm}

\par

For the proof we need the following lemma.

\par

\begin{proof}[Proof of Theorem \ref{fourop3}]
First assume that $p=1$, and
let $\chi \in C_0^\infty (\rr n)$ be such that $\int \chi _1\, dx=1$,
$\chi =\chi _2=\chi _1\otimes \chi _1$, $\chi _3=\chi _1\otimes \chi
_1\otimes \chi _1$ and $\psi \in C_0^\infty (\rr {3n})$ be such that
$\psi =1$ in $\supp \chi _3$. We also let $X_1=(x_1,y_1,\zeta _1)$,
$X=(x,y,\zeta )$, and consider the modulus of short-time Fourier
transform of the distribution kernel of $K_{a,\fy}$, i.{\,}e.
\begin{align*}
I_a(x,y,\xi ,\eta ) &= |\mathscr F(K_{a,\fy}\, \chi _2(\cdo -(x,y
))(\xi ,\eta )|
\\[1ex]
&= \Big |\int a(x_1,\zeta _1)e^{i\varphi ( X_1)}
\chi _2(x_1-x,y_1-y)e^{-i(\scal {x_1}\xi +\scal {y_1}\eta )}\, d X_1
\Big |,
\\[1ex]
&= \Big |\iint a(x_1,\zeta _1)e^{i\varphi ( X_1)}
\chi _3(X_1-X)e^{-i(\scal {x_1}\xi +\scal {y_1}\eta )}\, d X_1d\zeta
\Big |,
\end{align*}
and note that the $L^1_{(\omega _0)}$-norm of $I_a$ is equivalent to
the $M^1_{(\omega _0)}$-norm of $K_{a,\fy}$ in view of Remark
\ref{p1.7}. By a change of variables it follows that
\begin{equation*}
I_a(x,y,\xi ,\eta )=
\Big |\iint a(x_1+x,\zeta _1+\zeta )e^{i\varphi ( X_1+X)}
\chi _3(X_1)e^{-i(\scal {x_1}\xi +\scal {y_1}\eta )}\, d X_1d\zeta
\Big |
\end{equation*}
In a similar way as in Section \ref{sec2} we let $\psi _{1,X}$ and $\psi _{2,X}$ be defined as in \eqref{psidef}.
By letting $a_1(x,y,\xi) =a(x,\xi )$, an application of Taylor formula
on $\varphi$ gives
\begin{multline*}
I_a(x,y, \xi ,\eta )
\\[1ex]
= \Big |\iint a(x_1+x,\zeta _1+\zeta )e^{i\psi
_{2,X}(X_1)}\chi _3(X_1)e^{-i(\scal {x_1}\xi
+\scal {y_1}\eta -\psi _{1,X}(X_1))}\, d X_1d\zeta \Big |
\\[1ex]
= \Big |\int  e^{i\fy (X)}\mathscr F( a(\cdo +x,\cdo +\zeta )\, \chi
_3\, e^{i\psi _{2,X}})(\xi -\fy '_1(X),\eta -\fy '_2(X), -\fy
'_3(X))\,  d\zeta \Big |
\\[1ex]
\le \int ( |\mathscr F( a_1\, \chi _3(\cdo -X
))|*|\mathscr F(e^{i\psi _{2,X}})| (\xi -\fy '_1(X),\eta -\fy '_2(X),
-\fy '_3(X))\,  d\zeta
\end{multline*}
This implies that
\begin{equation}\label{decomp1}
I_a(x,y, \xi ,\eta ) \le \sum _{k=0}^\infty I_{a,k}(x,y,\xi ,\eta
)/k!,
\end{equation}
where
\begin{align*}
I_{a,0}(x,y,\xi ,\eta) &\equiv \int ( |\mathscr F( a_1\, \chi _3(\cdo
-X ))| (\xi -\fy '_1(X),\eta -\fy '_2(X), -\fy '_3(X))\,  d\zeta
\\[1ex]
I_{a,k}(x,y,\xi ,\eta) &\equiv \int ( |\mathscr F( a_1\, \chi _3(\cdo
-X ))|*\Phi _{k, X} (\xi -\fy '_1(X),\eta -\fy '_2(X),
-\fy '_3(X))\,  d\zeta ,
\\[1ex]
\Phi _{k, X} &\equiv |\mathscr F({\psi}_{2, X})| \ast
\cdots \ast |\mathscr F({\psi}_{2, X})|, \qquad  k\ge 1.
\end{align*}
Here the number of factors in the latter convolutions is equal to $k$.

\par

We need to estimate the $L^1_{(\omega _0)}$ norm of $I_{a,k}(x,y,\xi
,\eta)$, and start to consider the case $k=0$.

\vspace{1cm}


\par

Next we consider $I_{a,k}(x,y,\xi ,\eta)$ when $k\ge 1$. An application
of Lemma \ref{fourop3} shows that there is a function $G$ such
that $|\mathscr F({\psi}_{2, X})|\le G$ and $\nm G{L^1_{(v)}}\le
C\nm {\fy ''}{M^{\infty ,1}_{(v)}}$ for some constant $C>0$. Hence if
$\Upsilon _k\equiv G*\cdots *G$ with $k$ factors of $G$ in the
convolution, then it follows that $|\Phi _{k, X}|\le \Upsilon
_k$ and that $\nm {\Upsilon _k}{L^1_{(v)}}\le C^k\nm {\fy ''}{M^{\infty
,1}_{(v)}}^k$, where the latter inequality follows from the fact that
$v$ is submultiplicative.

\par

By letting $\xi _0=\xi -\fy '_1(X)$, $y_0=-\fy '_3(X)$ and
$$
J_k(x,y,\xi ,\eta ,z) = |\mathscr F( a_1\, \chi _3(\cdo
-X ))|*\Upsilon _{k} (\xi _0,\eta -\fy '_2(X), y_0),
$$
we therefore get
$$
I_{a,k}(x,y,\xi ,\eta )\le \int J_k(x,y,\xi ,\eta ,\zeta )\, d\zeta
$$
and
\begin{multline*}
J_{k}(x, y, \xi ,\eta ,\zeta )
\\[1ex]
\le \iiint |V_\chi a(x,\zeta ,\xi _0-\xi _1,y_0-y_1)|
|\widehat \chi (\eta -\fy '_2(X)-\eta _1)|\Upsilon _k(\xi _1,\eta
_1,y_1)\, d\xi _1 d\eta _1dy_1.
\end{multline*}
By \eqref{weightcond2} we have for some constant $C$ that
\begin{align*}
\omega _0(x,y,\xi ,\eta ) &\le C\omega _0(x,y,\xi ,\fy '_2(X))v_1(\eta
-\fy _2(X)-\eta _1)v_1(\eta _1)
\\[1ex]
\omega _0(x,y,\xi ,\fy _2(X)) &= \omega (x,\zeta ,\xi _0,y_0)
\\[1ex]
&\le C\omega (x,\zeta ,\xi _0-\xi _1,y_0-y_1)v_2(\xi _1,y_1),
\end{align*}
which implies that
$$
\omega _0(x,y,\xi ,\eta ) \le C^2\omega (x,\zeta ,\xi _0-\xi
_1,y_0-y_1)v_1(\eta -\fy _2(X)-\eta _1)v(\xi _1,\eta _1,y_1)
$$

\par

Hence if
\begin{align*}
F(x,\zeta ,\xi ,z) &= |V_\chi a(x,\zeta ,\xi ,z)\omega (x,\zeta ,\xi
,z)|
\\[1ex]
G(\eta ) &= |v_1(\eta )\widehat \chi (\eta )|
\\[1ex]
H_k(\eta ) &= \Upsilon _k(\xi ,\eta ,y)v_1(\eta )v_2(\xi ,y) =
\Upsilon _k(\xi ,\eta ,y)v(\xi ,\eta ,y),
\end{align*}
then we get
\begin{multline*}
J_k(x,y,\xi ,\eta ,\zeta )\omega _0(x,y,\xi ,\eta )
\\[1ex]
\le C\iiint F(x,\zeta ,\xi _0-\xi _1,y_0-y_1)G(\eta -\fy '_2(X)-\eta
_1)H_k(\xi _1,\eta _1,y_1)\, d\xi _1d\eta _1dy_1.
\end{multline*}

\par

By applying the $L^1$-norm on the latter estimate we get
\begin{multline*}
\nm {I_{a,k}}{L^1_{(\omega _0)}} \le C_1\int \nm {J_k(\cdo ,\zeta
)\omega _0}{L^1}\, d\zeta 
\\[1ex]
\le C_2\nm G{L^1}\nm {H_k}{L^1}\iiiint F(x,\zeta ,\xi ,-\fy '_3(X))\,
dxdyd\xi d\zeta 
\\[1ex]
\le C_3\nm {\fy ''}{M^{\infty ,1}_{(v)}}^k \iiiint F(x,\zeta ,\xi
,-\fy '_3(X))\, dxdyd\xi d\zeta ,
\end{multline*}
for some constants $C_1,\dots ,C_3$. Hence, by taking $(x,-\fy
'_3(X),\xi ,\zeta )$ as new variables of integration, and using
\eqref{phasecond} we get
$$
\nm {I_{a,k}}{L^1_{(\omega _0)}} \le C\mathsf d ^{-1}\nm F{L^1}\nm {\fy
''}{M^{\infty ,1}_{(v)}}^k.
$$

Hence, by applying the $L^1_{(\omega _0)}$ norm on the latter estimate
we get
\begin{multline*}
\nm {K_{a,\fy}}{M^{1}_{(\omega _0)}}\le \nm {I_a}{L^1}
\\[1ex]
=\sum _{k=0}^\infty \frac{1}{k!} \iiint (|\mathscr{F}(a \otimes \chi
)|\ast |\Upsilon _{k}|)(\xi -\partial _x\fy ( X_0),-\partial _\zeta
\fy ( X_0),\eta -\partial _y\fy
( X_0))\, dyd\xi d\eta
\\[1ex]
=\sum _{k=0}^\infty \frac{1}{k!} \iiint (|\mathscr{F}(a \otimes \chi
)|\ast |\Upsilon _{k}|)(\xi ,-\partial _\zeta
\fy ( X_0),\eta )\, dyd\xi
d\eta 
\\[1ex]
\le \mathsf d ^{-1}\sum _{k=0}^\infty \frac{1}{k!} \iiint
(|\mathscr{F}(a \otimes \chi )|\ast |\Upsilon _{k}|)(\xi ,x,\eta )\,
dxd\xi d\eta
\\[1ex]
\le C\mathsf d ^{-1}\sum _{k=0}^\infty \frac{1}{k!} \nm a{M^{1}}(C\nm
{\fy ''}{M^{\infty ,1}})^k = C\mathsf d ^{-1}\nm a{M^{1}}\exp (C\nm
{\fy ''}{M^{\infty ,1}}),
\end{multline*}
where the second inequality follows from \eqref{phasecond} and
taking $x=-\partial _\zeta \fy ( X_0)$ as new variable of
integration in the $y$-direction. This proves the assertion in this
case.

\par

For general $a\in M^{1}$, the asserted continuity now follows by
applying Proposition \ref{modspchar} in a way similar as in the proof
of Proposition \ref{fourop}. We leave the details for the reader.

\par

Next we consider the case $p=\infty$. Assume that $a\in M^1_{(\omega )}(\rr
{2n})$ and $b\in M^1_{(1/\omega )}(\rr
{2n})$, and let $\widetilde \fy (x,y,\xi )=-\fy (x,\xi ,y)$. Then
\eqref{phasecond} also holds when $\fy$ is replaced by $\widetilde
\fy$. Hence, the first part of the proof shows that $K_{b,\widetilde
\fy}\in M^{1}_{(1/\omega _0)}$. Furthermore, by straight-forward
computations we have
\begin{equation}\label{Kadjoint}
(K_{a,\fy},b)=(a,K_{b,\widetilde \fy}).
\end{equation}
In view of Proposition \ref{p1.4} (3), it follows that the right-hand
side in \eqref{Kadjoint} makes sense if, more generally, $a$ is an
arbitrary element in $M^{\infty}_{(\omega )}(\rr {2n})$, and then
$$
|(a,K_{b,\widetilde \fy})|\le C\mathsf d^{-1}\nm a{M^\infty _{(\omega
)}}\nm b{M^1_{(1/\omega )}}\exp (C\nm {\fy ''}{M^{\infty,1}_{(v)}}),
$$
for some constant $C$ which is independent of $\mathsf d$, $a\in
M^\infty$ and $b\in M^1$.

\par

Hence, by letting $K_{a,\fy}$ be defined as
\eqref{Kadjoint} when $a\in M^\infty$, it follows that $a\mapsto
K_{a,\fy}$ on $M^1$ extends to a continuous map on
$M^\infty$. Furthermore, since $\mathscr S$ is dense in $M^\infty$
with respect to the weak$^*$ topology, it follows that this extension
is unique. We have therefore proved the theorem for $p\in \{ 1,\infty
\}$.

\par

For general $p\in [1,\infty]$, the result now follows by
interpolation, using Theorem 4.1.2 in \cite{BL} and Proposition
\ref{interpolmod}.
\end{proof}

\par

\begin{equation}\tag*{(\ref{phasecond})$'$}
 s_2t_1-s_1t_2=1,\qquad
|\det(\fy ''_{y,\xi }(x,y,\xi ))|\ge \mathsf d
\end{equation}
\begin{equation}\tag*{(\ref{weightcond1})$'$}
\begin{aligned}
&\omega
_0(s_2x-t_2y,-s_1x+t_1y,t_1\xi+s_1\varphi '_y(x,y,\eta
),t_2\xi+s_2\varphi '_y(x,y,\eta ))
\\[1ex]
&= \omega (x,\eta ,\xi -\varphi '_x(x,y,\eta ),-\varphi
'_\eta(x,y,\eta ))
\end{aligned}
\end{equation}
and
\begin{equation}\tag*{(\ref{weightcond2})$'$}
\begin{alignedat}{2}
\omega _0(x,y,\xi+s_1\zeta , \eta+s_2\zeta ) &\le  \omega _0(x,y,\xi
,\eta )v_1(\zeta),&\quad x,y,\xi,\eta ,\zeta &\in \rr n
\\[1ex]
\omega (x,\eta ,\xi _1+\xi _2,y_1+y_2) &\le \omega (x,\eta ,\xi
_1,y_1)v_2(\xi _2,y_2),&\quad x,y_1,y_2, \xi _1,\xi _2,\eta &\in \rr n
\\[1ex]
v(x,y,\zeta,\xi ,\eta ,z) &= v_1(\eta )v_2(\xi ,z),&\quad x,y,z,\xi
,\eta ,\zeta &\in \rr n.
\end{alignedat}
\end{equation}

\par

\renewcommand{\rubrik}{Theorem \ref{fourop3}$'$}

\par

\begin{tom}
Assume that $p\in [1,\infty ]$, $s_j,t_j\in \mathbf R$ for
$j=1,2$, $\mathsf d >0$, $v\in \mathscr P(\rr {6n})$ is
submultiplicative, $\omega _0,\omega \in \mathscr P(\rr {4n})$ and
that $\varphi \in C(\mathbf {R}^{3n})$ are such that $\fy$ is
real-valued, $\partial ^{\alpha } \varphi \in M ^{\infty,1}_{(v)}$ for
$|\alpha | =2$ and \eqref{phasecond}$'$--\eqref{weightcond2}$'$ are
fulfilled. Then the map
$$
a\mapsto K_{a,\fy}(x,y)\equiv \int a(t_1x+t_2y,\xi
)e^{i\varphi (t_1x+t_2y,s_1x+s_2y,\xi )}\, d\xi ,
$$
from $\mathscr S(\rr {2n})$ to $\mathscr S'(\rr {2n})$ extends
uniquely to a continuous map from $M^p_{(\omega )}(\rr {2n})$ to
$M^p_{(\omega _0)}(\rr {2n})$.
\end{tom}

\par

\begin{proof}
By letting
$$
x_1 =t_1x+t_2y,\quad y_1= s_1x+s_2y
$$
as new coordinates, it follows that we may assume that $t_1=s_2=1$ and
$t_2=s_1=0$, and then the result agrees with Theorem
\ref{fourop3}. The proof is complete.
\end{proof}

\par

Assume that $a\in M^{\infty}(\rr {2n})$, $t_1,t_2\in \mathbf R$, and
that $\fy \in C(\rr {3n})$ is real-valued and satisfies $\partial
^{\alpha } \varphi \in M ^{\infty,1}$ for $|\alpha | =2$ and
\eqref{phasecond} for some $\mathsf d >0$. Then we let the Fourier
integral operator $\opp _{\fy}(a)=\opp _{\fy ,t_1,t_2}(a)$ be the
continuous operator from $\mathscr S(\rr n)$ to $\mathscr S'(\rr n)$
with kernel $K_{a,\fy}$ in Theorem \ref{fourop3}. Furthermore, since
the case $t_1=1$ and $t_2=0$ is especially important we set $\opp
_{\fy ,0}(a)=\opp _{\fy ,1,0}(a)$. The following result is now an
immediate consequence of Theorem \ref{fourop3} and Theorem 4.3 in
\cite{To7}.

\par

Next we consider Fourier operators when $a\in M^{\infty ,1}_{(\omega
)}(\rr {2n})$. In the following it is natural to consider weights
$\omega _1,\omega _2\in \mathscr P(\rr {2n})$ and $\omega \in \mathscr
P(\rr {4n})$ which satisfy
\begin{equation}\label{weigthsineq1}
\frac {\omega _2(s_2x-t_2y,t_1\xi +s_1\fy '_y(x,y,\eta ))}{\omega
_1(-s_1x+t_1y,-t_2\xi -\fy '_2(x,y,\eta ))}\le C\omega (x,\eta ,\xi
-\fy '_1(x,y,\eta ),-\fy '_3(x,y,\eta )),
\end{equation}
for some constant $C$.

\par

The proof of the following proposition in the case $p=\infty$ can be
found in \cite{Gc2, To7A}.

\par

\begin{prop}\label{fourop3A}
Assume that  $p\in [1,\infty]$, $\omega_j\in \mathscr{P}(\rr {2n_j})$,
for $j=1,2$, and $\omega\in \mathscr{P}(\rr {2n_1+2n_2})$ fulfill for
some positive constant $C$
$$
\frac{\omega_2(x,\xi)}{\omega_1(y,-\eta)}\le C \omega(x,y,\xi,\eta). 
$$
Assume moreover that  $K\in M^p_{(\omega)}(\rr {n_1+n_2})$ and $T$ is
the linear and continuous  map from $\mathscr S(\rr {n_1})$ to
$\mathscr S'(\rr {n_2})$ defined by:
$$
(Tf)(x)= (K(x,\cdot), f), \quad f\in \mathscr S(\rr {n_1}).
$$
Then $T$ extends uniquely to a continuous  map from
$M^{p'}_{(\omega_1)}(\rr{n_1})$ to $M^p_{(\omega_2)}(\rr{n_2})$
\end{prop}

\par

\begin{proof}
By Proposition \ref{p1.4} (3) and duality, it sufficies to prove that
for some constant $C$ independent of  $f\in \mathscr S(\rr {n_1})$ and
$g\in \mathscr S(\rr{n_2})$, it holds:
$$
|(K,g\otimes \bar f)|\le C\Vert K\Vert_{M^p_{(\omega)}}
\Vert g\Vert_{M^{p'}_{(1/\omega_2)}} \Vert f
\Vert_{M^{p'}_{(\omega_1)}}.
$$
Let us set $\omega_3(x,\xi)=\omega_1(x,-\xi)$, then by straightforward
calculation and using Remark \ref{p1.7} (7) we get
$$
\begin{array}{ll}
|(K,g\otimes \bar f)|&\le C_1\Vert K  \Vert_{M^{p}_{(\omega)}} \Vert
g\otimes \bar f \Vert_{M^{p'}_{(1/\omega)}} \le C_2 \Vert K
\Vert_{M^{p}_{(\omega)}} \Vert g  \Vert_{M^{p'}_{(1/\omega_2)}} \Vert
\bar f \Vert_{M^{p'}_{(\omega_3)}}
\\[1ex]
& \le C \Vert K  \Vert_{M^{p}_{(\omega)}}
\Vert g  \Vert_{M^{p'}_{(1/\omega_2)}} \Vert  f
\Vert_{M^{p'}_{(\omega_1)}}
\end{array}
$$
\end{proof}

\par

\begin{prop}\label{fourop4}
Assume that,  $1<p<\infty$,     
$\mathsf d >0$, $v\in \mathscr P(\rr {6n})$ is
submultiplicative, $\omega _1,\omega _2\in \mathscr P(\rr {2n})$ and
$\omega \in \mathscr P(\rr {4n+2m})$ are such that \eqref{detphicond}
is fulfilled for some constant $C$, and that $a\in M^{\infty
,1}_{(\omega )} (\rr {2n+m})$. Also assume that $\varphi \in C(\mathbf
{R}^{3n})$ are such that $\fy$ is real-valued, $\partial ^{\alpha }
\varphi \in M ^{\infty,1}_{(v)}$ for $|\alpha | =2$ and
\eqref{detphicond} are fulfilled. Then $\opp _{\fy ,0}(a)$ extends to
a continuous operator from $M^p_{(\omega _1)}(\rr n)$ to $M^p_{(\omega
_2)}(\rr n)$, and
$$
\nm {\opp _{\fy ,0}(a)f}{M^p_{(\omega _2)}}\le C\mathsf d ^{-1}\nm
a{M^{\infty ,1}_{(\omega )}}\nm f{M^{p}_{(\omega _1)}}\exp (C\nm {\fy
''}{M^{\infty 1,}_{(v)}}).
$$
\end{prop}

\par

\begin{proof}
We shall mainly follow the proof of Theorem 3.2?? in \cite{Bo}. Assume
that $f,g\in \mathscr S(\rr n)$, and that $0\le \chi ,\psi \in
C_0^\infty (\rr {2n+m})$ and $0\le \chi _0\in C_0^\infty (\rr n)$ are
such that
$$
\nm {\chi _0} {L^1}=\nm {\chi }{L^2}=1
$$
and $\psi =1$ on $\supp \chi $. Also let $X=(x,y,\zeta )\in \rr {2n+m}$
and $X_1=(x_1,y_1,\zeta _1)\in \rr {2n+m}$ as usual. By
straight-forward computations we get
\begin{multline*}
(\op _\fy (a)f,g) = \int a(X)f(y)\overline {g(x)} e^{i\fy (X )}\, dX
\\[1ex]
= \iint a(X+X_1)\chi (X_1)^2f(y+y_1)\chi _0(y_1)\overline
{g(x+x_1)\chi _0(x_1)} e^{i\psi (X_1)\fy (X+X_1 )}\, dXdX_1
\end{multline*}
If $\mathscr F _{1,2}a$ denotes the partial Fourier transform of
$a(x,y,\zeta )$ with respect to the $x$ and $y$ variables, then
Parseval's formula gives
\begin{multline*}
(\opp _\fy (a)f,g)
\\[1ex]
=\iiiint  F(X,\xi ,\eta ,\zeta _1) \mathscr F(f(y+\cdo )\chi _0)(-\eta
)\overline {\mathscr F(g(x+\cdo )\chi _0)(\xi )}\, dXd\xi d\eta d\zeta
_1
\\[1ex]
=\iiiint  F(X,\xi ,\eta ,\zeta _1) (V_{\chi _0}f)(y,-\eta )\overline
{(V_{\chi _0}g)(x,\xi )}e^{-i(\scal x\xi+\scal y{\eta})}\, dXd\xi
d\eta d\zeta _1,
\\[1ex]
=\iiint \Big ( \int  F(X,\xi ,\eta ,\zeta _1)\, d\zeta _1 \Big )
(V_{\chi _0}f)(y,-\eta )\overline {(V_{\chi _0}g)(x,\xi )}e^{-i(\scal
x\xi+\scal y{\eta})}\, dXd\xi d\eta ,
\end{multline*}
where
\begin{equation*}
F(X,\xi ,\eta ,\zeta _1) = \mathscr F_{1,2}\big (e^{i\psi (\cdo ,\zeta
_1)\fy (X+(\cdo ,\zeta _1 ))}a(X+(\cdo ,\zeta _1))\chi (\cdo ,\zeta
_1)^2\big )(\xi ,\eta ).
\end{equation*}

\par

In order to further reformulate the action of $\opp _\fy (a)$, we let
$\psi _{1,X}$ and $\psi _{2,X}$ be the same as in \eqref{psidef}, and
set
\begin{equation}\label{Hdef}
\begin{aligned}
&\mathcal H_{a,\fy}(X,\xi ,\eta )
\\[1ex]
&= h_X*(\mathscr F(a(\cdo +X)\chi ))(\xi -\fy '_x(X),\eta -\fy
'_y(X),-\fy '_\zeta ),\quad \text{where}
\\[1ex]
h_X &= (2\pi )^{-(n+m/2)} (\mathscr F(e^{i\psi _{2,X}}\chi ))
\end{aligned}
\end{equation}

Next let $\psi _{1,X}$ and $\psi _{2,X}$ be the same as in
\eqref{psidef}. Then
\begin{multline*}
\int F(X,\xi ,\eta ,\zeta _1)\, d\zeta _1
\\[1ex]
=\mathscr F((e^{i\psi _{2,X}}\chi ) (a(\cdo +X)\chi )(\xi -\fy
'_x(X),\eta -\fy '_y(X),-\fy '_\zeta (X))
\\[1ex]
=\mathcal H_{a,\fy}(X,\xi -\fy '_x(X),\eta -\fy
'_y(X),-\fy '_\zeta (X)),
\end{multline*}
where
\begin{equation}\label{Hdef2}
\begin{aligned}
&\mathcal H_{a,\fy}(X,\xi ,\eta )
\\[1ex]
&= (2\pi )^{-(n+m/2)} (\mathscr F(e^{i\psi _{2,X}}\chi ))*(\mathscr
F(a(\cdo +X)\chi ))(X,\xi -\fy '_x(X),\eta -\fy '_y(X),-\fy '_\zeta
(X)).
\end{aligned}
\end{equation}

\par

Summing up we have proved that
\begin{equation}\label{fourrel}
(\opp _\fy (a)f,g) = \iiint  \mathcal
H_{a,\fy}(X,\xi ,\eta )
(V_{\chi _0}f)(y,-\eta )\overline {(V_{\chi _0}g)(x,\xi )}e^{-i(\scal
x\xi+\scal y{\eta})}\, dXd\xi d\eta
\end{equation}
Since
$$
|\mathscr F(a(\cdo +X)\chi )(\xi ,\eta ,z)| = |(V_\chi a)(X,\xi ,\eta
,z)|
$$
we therefore get
\begin{multline*}
\Big | \int F(X,\xi ,\eta ,\zeta _1)\, d\zeta _1 \Big | \le \big
(|\mathscr F(e^{i\psi _{2,X}}\chi ))|*|V_\chi a(X,\cdo )|\big )(\xi
-\fy '_x(X),\eta -\fy '_y(X),-\fy '_\zeta (X))
\\[1ex]
\le \sum _{k=0}^\infty \frac 1{k!} \big (|\widehat \chi | *\Phi _{k,X}
* |V_\chi a(X,\cdo )|\big )(\xi -\fy '_x(X),\eta -\fy '_y(X),-\fy
'_\zeta (X)),
\end{multline*}
where
\begin{align*}
\Phi _{0,X} &= \delta _0
\\[1ex]
\Phi _{k,X} &= |\mathscr F(\psi _{2,X})| * \cdots *|\mathscr F(\psi
_{2,X})|,\qquad k\ge 1,
\end{align*}
with $k$ factors in the latter convolution.

\par

By Lemma \ref{fourop3} it follows that there is a function $G\in
L^1_{(v)}$ such that $|\mathscr F({\psi}_{2, X})|\le G$ and $\nm
G{L^1_{(v)}}\le C\nm {\fy ''}{M^{\infty ,1}_{(v)}}$ for some constant
$C>0$. Hence if
$\Upsilon _k\equiv G*\cdots *G$ with $k$ factors of $G$ in the
convolution, then it follows that $|\Phi _{k, X}|\le \Upsilon
_k$ and that
$$
\nm {\Upsilon _k}{L^1_{(v)}}\le C_1^k\nm {G}{L^1_{(v)}}^k
\le C_2^k\nm {G}{M^1_{(v)}}^k \le C_3^k\nm {\fy ''}{M^{\infty
,1}_{(v)}}^k,
$$
for some constants $C_1,\dots ,C_3$. Here the first inequality in the
latter estimate follows from the fact that $v$ is submultiplicative.

\par

By letting
\begin{align*}
J_0(X,\xi ,\eta ) &= \big (|\widehat \chi |* |V_\chi a(X,\cdo )|\big
)(\xi -\fy '_x(X),\eta -\fy '_y(X),-\fy '_\zeta (X))
\\[1ex]
J_k(X,\xi ,\eta ) &= \big (|\widehat \chi |* \Upsilon _k* |V_\chi
a(X,\cdo )|\big )(\xi -\fy '_x(X),\eta -\fy '_y(X),-\fy '_\zeta
(X)),\quad k\ge 1,
\end{align*}
it follows now that
$$
\Big | \int F(X,\xi ,\eta ,\zeta _1)\, d\zeta _1 \Big | \le C \sum
_{k=0}^\infty J_k(X,\xi ,\eta ).
$$

Let $a_1(x,y,\zeta )=a(x,\zeta )$. For some real-valued function
$\Theta $, we have
\begin{multline*}
|(\opp _{\fy ,0}(a)f,g)| = (2\pi )^{-n}\Big | \iiint e^{i\fy }V_{\chi _3}a_1(X,\xi ,\eta ,-\fy '_3(X))\dot {}
\\[1ex]
V_\chi f(y,-\eta -\fy '_2(X))\overline{V_\chi g(x,\xi +\fy '_1(X))}e^{i\Theta (X,\xi ,\eta )}\, d\xi d\eta dX\Big |
\\[1ex]
\le C\iiint | V_{\chi _2}a(x,\zeta ,\xi ,-\fy '_3(X)) \, \widehat \chi (\eta ) V_\chi f(y,-\eta -\fy '_2(X))\overline{V_\chi g(x,\xi +\fy '_1(X))}|\, d\xi d\eta dX
\end{multline*}
In order to include the conditions for weight functions we set
\begin{align*}
H(\xi ,y) &= \sup _{x,\zeta }|(V_{\chi _2}a)(x,\zeta ,\xi ,y)\omega (x,\zeta ,\xi ,y)|
\\[1ex]
F_1(x,\xi ) &= |(V_\chi f)(x,\xi )\omega _1(x,\xi )| ,\quad F_2(x,\xi ) = |(V_\chi g)(x,\xi )(\omega _2(x,\xi ))^{-1}|,
\end{align*}
and we note that
$$
\nm H = \nm a{M^{\infty ,1}_{(\omega )}},\quad \nm {F_1}{L^p}=\nm f{M^p_{(\omega _1)}},\quad \nm {F_1}{L^{p'}}=\nm f{M^{p'}_{(1/\omega _2)}}.
$$
By taking  $\xi +\fy '_1(X)$, $\eta$ and $X$ as new variables of integration, and using \eqref{weigthsineq1} we obtain
\begin{multline*}
|(\opp _{\fy ,0}(a)f,g)| \le C \iiint | V_{\chi _2}a(x,\zeta ,\xi ,-\fy '_3(X)) \omega (x,\zeta ,\xi -\fy '_1(X),-\fy '_3(X))|\cdot {}
\\[1ex]
|V_\chi f(y,-\eta -\fy '_2(X))\omega _1(y,-\fy '_2(X))| |V_\chi g(x,\xi +\fy '_1(X))(\omega _2(x,\xi ))^{-1}| |\widehat \chi (\eta )|\, d\xi d\eta dX
\\[1ex]
\le C \iiint H(\xi -\fy '_1(X),-\fy '_3(X))F_2(x,\xi )| V_\chi f(y,-\eta -\fy '_2(X))\omega _1(y,-\fy '_2(X)) |\, |\widehat \chi (\eta )|\, d\xi d\eta dX.
\end{multline*}
Since $\omega _1$ belongs to $\mathscr P(\rr {2n})$, it follows that
$$
\omega _1(y,-\fy '_2(X)) \le \omega _1(,y,-\eta -\fy '_2(X))v_1(\eta ),
$$
for some $v_1\in \mathscr P(\rr n)$, giving that
\begin{multline*}
V_\chi f(y,-\eta -\fy '_2(X))\omega _1(y,-\fy '_2(X)) |\, |\widehat \chi (\eta )|
\\[1ex]
\le V_\chi f(y,-\eta -\fy '_2(X))\omega _1(y,-\eta -\fy '_2(X)) |\, |\widehat \chi (\eta )v_1(\eta )|\le F_1(y,-\eta -\fy '_2(X))h(\eta ),
\end{multline*}
where $h(\eta )=|\widehat \chi (\eta )v_1(\eta )|\in L^1$.

\par

A combination of these estimates and H{\"o}lder's inequality give
\begin{equation}\label{eq11}
|(\opp _{\fy ,0}(a)f,g)| \le C\iiint H(\xi -\fy '_1(X),-\fy '_3(X))F_1(y,-\eta - \fy '_2(X))F_2(x,\xi )h(\eta )\, d\xi d\eta dX \le CJ_1\cdot J_2,
\end{equation}
where
\begin{align*}
J_1^p &= \iiint H(\xi -\fy '_1(X),-\fy '_3(X))F_1(y,-\eta - \fy '_2(X))^ph(\eta )\, d\xi d\eta dX,
\\[1ex]
J_2^{p'} &= \iiint H(\xi -\fy '_1(X),-\fy '_3(X)) F_2(x,-\xi)^{p'}h(\eta )\, d\xi d\eta dX
\end{align*}

\par

We have to estimate $J_1$ and $J_2$. By taking $x_0 =\fy '_3(X)$,
$\zeta _0=\fy '_2(X)$, $y$, $\xi$ and $\eta$ as new variables of
integrations, and using \eqref{detphicond}, it follows that
\begin{multline*}
J_1 \le \Big ( \mathsf d ^{-1}\int \cdots \int  H(\xi -\fy
'_1(X),-x_0) F_1(y,-\eta - \zeta _0)^p h(\eta )\, dx_0dyd\xi d\eta
d\zeta _0 \Big )^{1/p}
\\[1ex]
=(\mathsf d ^{-1}\nm H{L^1}\nm h{L^1})^{1/p}\nm {F_1}{L^p}.
\end{multline*}
Hence
\begin{equation}\label{eq22}
J_1\le  (C_h\mathsf d ^{-1}\nm a{M^{\infty ,1}_{(\omega )}})^{1/p}\nm
{f}{M^p_{(\omega _1)}} .
\end{equation}
If we instead take $x$, $y_0=\fy '_3(X)$, $\xi $, $\eta $ and $\zeta
_0=\fy '_1(X)$ as new variables of integrations, it follows by similar
arguments that
\begin{equation}\tag*{(\ref{eq22})$'$}
J_2\le (C_h\mathsf d ^{-1}\nm a{M^{\infty ,1}_{(\omega )}})^{1/p'}\nm
{g}{M^{p'}_{(1/\omega _2)}} .
\end{equation}
A combination of \eqref{eq11}, \eqref{eq22} and \eqref{eq22}$'$ now give
$$
|(\opp _{\fy ,0}(a)f,g)| \le  C_h\mathsf d ^{-1}\nm a{M^{\infty ,1}_{(\omega )}} \nm {f}{M^p_{(\omega _1)}} \nm {g}{M^{p'}_{(1/\omega _2)}},
$$
which proves the assertion.
\end{proof}

\par

\begin{thm}\label{thm4}
Assume that $p\in [1,\infty ]$, $a\in M^{p}(\rr {2n})$,
$t_1,t_2\in \mathbf R$, and that $\fy \in C(\rr {3n})$ is real-valued
and satisfies $\partial ^{\alpha } \varphi \in M ^{\infty,1}$ for
$|\alpha | =2$ and \eqref{phasecond} for some $\mathsf d >0$. Then the
definition of $\opp _{\fy ,t_1,t_2}(a)$ from $\mathscr S(\rr n)$ to
$\mathscr S'(\rr n)$ extends uniquely to a continuous map from
$M^{p'}(\rr n)$ to $M^p(\rr n)$.
\end{thm}

\par

By combining Theorem \ref{fourop2}, Theorem \ref{fourop3} and
interpolation, we obtain the following result.

\par

\begin{thm}\label{thm5}
Assume that $p,q\in [1,\infty ]$ are such that $q\le \min (p,p')$,
$a\in M^{p,q}(\rr {2n})$, $t_1,t_2\in \mathbf R$, and that $\fy \in
C(\rr {3n})$ is real-valued and satisfies $\partial ^{\alpha } \varphi
\in M ^{\infty,1}$ for $|\alpha | =2$. Also assume that
\eqref{detphicond} and \eqref{phasecond} are fulfilled for some
$\mathsf d >0$. Then the definition of $\opp _{\fy ,t_1,t_2}(a)$ from
$\mathscr S(\rr n)$ to $\mathscr S'(\rr n)$ extends uniquely to a
continuous map from $M^{p'}(\rr n)$ to $M^p(\rr n)$. Then the
definition of $\opp _{\fy ,t_1,t_2}(a)$ from $\mathscr S(\rr n)$ to
$\mathscr S'(\rr n)$ extends uniquely to a Schatten-von Neumann
operator of order $p$ on $L^2(\rr n)$.
\end{thm}

\par

\begin{proof}
We may assume that $q=\min (p,p')$. First assume that $p\le 2$, and
let $b\in \mathscr S'(\rr {2n})$ be chosen such that $b(x,D)=\opp
_{\fy}(a)$. Then the operator kernel of $b$ belongs to $M^p$, and
since $M^p$ is invariant under partial Fourier transformations in view
of Remark \ref{p1.7} (5), the result is a consequence of Proposition
1.7 in \cite{To7}.

\par

If instead $p=\infty$, then it follows from Theorem \ref{fourop2} that
$\opp _{\fy}(a)$ is continuous on $L^2$, which proves the result in
this case as well. The result now follows for general $p\in [2,\infty
]$ by interpolation, using Proposition \ref{interpolmod} and
\eqref{interpschatt}. The proof is complete.
\end{proof}

\end{document}